\documentclass[12pt]{article}
\usepackage{
latexsym, amssymb, amscd, amsthm, graphicx, 
}

 \title{{\bf Meromorphic open-string vertex algebras}}
 \author{Yi-Zhi Huang}
    \date{}
    \begin{document}
    \bibliographystyle{alpha}
    \maketitle
\newtheorem{thm}{Theorem}[section]
\newtheorem{defn}[thm]{Definition}
\newtheorem{prop}[thm]{Proposition}
\newtheorem{cor}[thm]{Corollary}
\newtheorem{lemma}[thm]{Lemma}
\newtheorem{rema}[thm]{Remark}
\newtheorem{app}[thm]{Application}
\newtheorem{prob}[thm]{Problem}
\newtheorem{conv}[thm]{Convention}
\newtheorem{conj}[thm]{Conjecture}
\newtheorem{cond}[thm]{Condition}
    \newtheorem{exam}[thm]{Example}
\newtheorem{assum}[thm]{Assumption}
     \newtheorem{nota}[thm]{Notation}
\newcommand{\halmos}{\rule{1ex}{1.4ex}}
\newcommand{\pfbox}{\hspace*{\fill}\mbox{$\halmos$}}

\newcommand{\nn}{\nonumber \\}

 \newcommand{\res}{\mbox{\rm Res}}
 \newcommand{\ord}{\mbox{\rm ord}}
\renewcommand{\hom}{\mbox{\rm Hom}}
\newcommand{\edo}{\mbox{\rm End}\ }
 \newcommand{\pf}{{\it Proof.}\hspace{2ex}}
 \newcommand{\epf}{\hspace*{\fill}\mbox{$\halmos$}}
 \newcommand{\epfv}{\hspace*{\fill}\mbox{$\halmos$}\vspace{1em}}
 \newcommand{\epfe}{\hspace{2em}\halmos}
\newcommand{\nord}{\mbox{\scriptsize ${\circ\atop\circ}$}}
\newcommand{\wt}{\mbox{\rm wt}\ }
\newcommand{\swt}{\mbox{\rm {\scriptsize wt}}\ }
\newcommand{\lwt}{\mbox{\rm wt}^{L}\;}
\newcommand{\rwt}{\mbox{\rm wt}^{R}\;}
\newcommand{\slwt}{\mbox{\rm {\scriptsize wt}}^{L}\,}
\newcommand{\srwt}{\mbox{\rm {\scriptsize wt}}^{R}\,}
\newcommand{\clr}{\mbox{\rm clr}\ }
\newcommand{\tr}{\mbox{\rm Tr}}
\newcommand{\C}{\mathbb{C}}
\newcommand{\Z}{\mathbb{Z}}
\newcommand{\R}{\mathbb{R}}
\newcommand{\Q}{\mathbb{Q}}
\newcommand{\N}{\mathbb{N}}
\newcommand{\CN}{\mathcal{N}}
\newcommand{\F}{\mathcal{F}}
\newcommand{\I}{\mathcal{I}}
\newcommand{\V}{\mathcal{V}}
\newcommand{\one}{\mathbf{1}}
\newcommand{\BY}{\mathbb{Y}}
\newcommand{\ds}{\displaystyle}

        \newcommand{\ba}{\begin{array}}
        \newcommand{\ea}{\end{array}}
        \newcommand{\be}{\begin{equation}}
        \newcommand{\ee}{\end{equation}}
        \newcommand{\bea}{\begin{eqnarray}}
        \newcommand{\eea}{\end{eqnarray}}
         \newcommand{\lbar}{\bigg\vert}
        \newcommand{\p}{\partial}
        \newcommand{\dps}{\displaystyle}
        \newcommand{\bra}{\langle}
        \newcommand{\ket}{\rangle}

        \newcommand{\ob}{{\rm ob}\,}
        \renewcommand{\hom}{{\rm Hom}}

\newcommand{\A}{\mathcal{A}}
\newcommand{\Y}{\mathcal{Y}}

\newcommand{\dlt}[3]{#1 ^{-1}\delta \bigg( \frac{#2 #3 }{#1 }\bigg) }

\newcommand{\dlti}[3]{#1 \delta \bigg( \frac{#2 #3 }{#1 ^{-1}}\bigg) }

\vspace{2em}



\renewcommand{\theequation}{\thesection.\arabic{equation}}
\renewcommand{\thethm}{\thesection.\arabic{thm}}
\setcounter{equation}{0} \setcounter{thm}{0}
\maketitle

\begin{abstract}
A notion of meromorphic open-string vertex algebra is introduced.
A meromorphic open-string vertex algebra is an open-string vertex 
algebra in the sense of Kong and the author satisfying additional 
rationality (or meromorphicity) conditions for vertex operators. 
The vertex operator map for a meromorphic open-string vertex algebra satisfies
rationality and 
associativity but in general does not satisfy 
the Jacobi identity, commutativity, the commutator formula, the skew-symmetry
or even the associator formula.
Given a vector space $\mathfrak{h}$, we construct a meromorphic 
open-string vertex algebra structure on the 
tensor algebra of the negative part of the affinization of $\mathfrak{h}$ such that 
the vertex algebra struture on the symmetric algebra
of the negative part of the Heisenberg algebra associated to
$\mathfrak{h}$ is a quotient of this meromorphic open-string vertex algebra.
We also introduce the notion of
left module for a meromorphic open-string vertex algebra
and construct left  modules for the meromorphic open-string vertex algebra
above. 
\end{abstract}

\tableofcontents

\section{Introduction}

Vertex (operator) algebras arose naturally in the study of two-dimensional 
conformal field theories in physics (see the first systematic study 
using the method of operator product expansion in \cite{BPZ} by Belavin, Polyakov and
Zamolodchikov) and in the vertex operator
construction of representations of affine Lie algebras and in the
construction and study of the ``moonshine module'' for the Monster
finite simple group in mathematics 
(see the announcement \cite{B} by Borcherds and 
the monograph \cite{FLM} by Frenkel, Lepowsky and Meurman \cite{FLM}). 

 Vertex (operator) algebras can be viewed as the
``closed-string-theoretic'' analogues of both Lie algebras and 
commutative associative algebras. A vertex (operator) algebra
is defined in terms of either the Jacobi identity or the duality property or 
parts of these axioms.
The Jacobi identity contains the commutator formula for vertex operators
and the duality property includes in particular commutativity.
The commutator formula and  commutativity are fundamental to vertex (operator) 
algebras. Many of the results on vertex (operator) algebras and their representations depend
heavily on the commutator formula and commutativity. 
The commutator formula and commutativity also play an important role in the construction 
of examples of vertex (operator) algebras, especially in the construction of 
vertex operator algebras associated to affine Lie algebras and the Virasoro algebra. 
It was proved in \cite{FHL}
that associativity and thus other properties, including in particular the Jacobi identity,
follows from commutativity and other minor axioms. Geometrically, it was shown 
in \cite{H0} and \cite{H1} by the author that commutativity is equivalent to 
a meromorphicity property on an algebra over the partial operad of the moduli
space of spheres with punctures and standard local coordinates at the punctures.

Beyond topological field theories, two-dimensional conformal field theories are
the only mathematically successful quantum field theories. Many people attribute this 
success to the existence of the infinite-dimensional conformal symmetry. But
from the experience in the study of two-dimensional conformal field theories 
in terms of the representation theory of vertex operator algebras, 
at least in the genus-zero case, commutativity 
or other equivalent properties is the main reason why two-dimensional conformal 
field theories are mathematically much better understood
than other non-topological quantum field theories. In fact, it is the commutator formula 
that allows one to apply the Lie-theoretic method to the study of vertex operator algebras
and their representations. 

Two-dimensional conformal field theories are deep mathematical theories
that play important roles in both mathematics and physics. But there are also 
other non-topological quantum field theories of fundamental importance
in both mathematics and physics. Nonlinear sigma models 
whose target spaces are not Calabi-Yau manifolds are such examples. 
The most important example
is Yang-Mills theory. In physics, Yang-Mills theory is known to 
describe fundamental interactions and in mathematics, one of the 
major unsolved problem is the 
existence of quantum Yang-Mills theory and the mass gap conjecture. 
Unfortunately, for these theories, we do not expect that they have a commutativity property 
as strong as the commutativity property 
for two-dimensional conformal field theories. 
This is actually the reason why it is very difficult to generalize the veretx-algebraic
approach in the study of two-dimensional conformal field theories to other non-topological 
quantum field theories. 

On the other hand, we have another fundamental property of vertex (operator) algebras:
Associativity. Associativity for vertex (operator) algebras 
is a strong form of operator product expansion
for meromorphic fields. While we do not expect that commutativity holds for 
general non-topological quantum field theories, we do believe that operator product 
expansion holds for these theories. Therefore, to study non-topological
quantum field theories that are not two-dimensional conformal field theories, 
one approach is to find algebraic structures satisfying certain associativity property
but not necessarily any commutativity property. 

In dimension $2$, Kong and the author introduced and constructed 
open-string vertex algebras
in \cite{HK}. An open-string vertex algebra satisfies associativity but not 
commutativity. However, the examples given in \cite{HK} are constructed
using modules and intertwining operators for a vertex operator algebra belonging to
the meromorphic
center of the open-string vertex algebra. Since intertwining operators
for the meromorphic center satisfy the commutativity property for intertwining operators
(a generalization of commutativity for vertex (operator) algebras formulated 
and proved in \cite{H6},
\cite{H3} and \cite{H5}), these examples 
of open-string vertex algebras still satisfy a certain generalized version of commutativity.
In fact, these open-string vertex algebra are still part of an open-closed two-dimensional 
conformal field theory describing the interaction of boundary states or open strings. 
To go beyond conformal field theories in dimension $2$, we need to find examples 
of open-string vertex algebras that are not constructed from modules and intertwining
operators for vertex (operator) algebras. 

In the present paper, we construct a class of such examples. In fact, the examples that
we construct in the present paper satisfy stronger conditions than those 
open-string vertex algebras constructed in \cite{HK}. Like vertex (operator) algebras,
the products and iterates of vertex
operators of these open-string vertex algebras are expansions of rational 
functions. We call an open-string vertex algebra satisfying such a rationality 
property a meromorphic open-string vertex algebra. The vertex operator 
map for a meromorphic open-string vertex algebra satisfies
rationality and 
associativity but in general does not satisfy 
the Jacobi identity, commutativity, the commutator formula, the skew-symmetry
or even the associator formula.

Given a 
vector space $\mathfrak{h}$, we have the Heisenberg algebra 
$\hat{\mathfrak{h}}=\hat{\mathfrak{h}}_{+}
\oplus \hat{\mathfrak{h}}_{-}\oplus \hat{\mathfrak{h}}_{0}$, where
$\hat{\mathfrak{h}}_{+}=\mathfrak{h}\otimes t\C[t]$, $\hat{\mathfrak{h}}_{-}=\mathfrak{h}
\otimes t^{-1}\C[t^{-1}]$,
$\hat{\mathfrak{h}}_{0}=\mathfrak{h}\oplus \C\mathbf{k}$  and 
$\C\mathbf{k}$ is the center of $\hat{\mathfrak{h}}$. 
Instead of the universal enveloping algebra 
$U(\hat{\mathfrak{h}})$ of $\hat{\mathfrak{h}}$, we consider 
the quotient $N(\hat{\mathfrak{h}})$ of the tensor algebra
$T(\hat{\mathfrak{h}})$ of $\hat{\mathfrak{h}}$ by only the commutator relations between 
$\hat{\mathfrak{h}}_{+}$ and $\hat{\mathfrak{h}}_{-}$, between 
$\hat{\mathfrak{h}}_{+}$ and $\hat{\mathfrak{h}}_{0}$, between 
$\hat{\mathfrak{h}}_{-}$ and $\hat{\mathfrak{h}}_{0}$ and between 
$\mathfrak{h}$ and $\C\mathbf{k}$, but not the commutator relations 
between $\hat{\mathfrak{h}}_{+}$ and itself, between $\hat{\mathfrak{h}}_{-}$ and itself 
and between $\mathfrak{h}$ and itself. In other words, we do not assume that 
$\mathfrak{h}$ and consequently $\hat{\mathfrak{h}}_{+}$ and $\hat{\mathfrak{h}}_{-}$
are abelian Lie algebras. Actually, since we work with 
the tensor algebras of $\hat{\mathfrak{h}}_{+}$ and $\hat{\mathfrak{h}}_{-}$
and $\mathfrak{h}$, we do not assume any relations among linearly independent 
elements of these vector spaces. From a left module $M$ for the tensor algebra 
$T(\mathfrak{h})$ of $\mathfrak{h}$, we 
construct an induced left module for $N(\hat{\mathfrak{h}})$ and prove that 
this induced left module is linearly isomorphic to $T(\hat{\mathfrak{h}}_{-})\otimes M$
where $T(\hat{\mathfrak{h}}_{-})$ is the tensor algebra of $\hat{\mathfrak{h}}_{-}$.
In the case that $M$ is the trivial left module $\C$ for $T(\mathfrak{h})$,
we construct a meromorphic open-string vertex algebra structure on 
$T(\hat{\mathfrak{h}}_{-})\simeq T(\hat{\mathfrak{h}}_{-})\otimes \C$. 
We know that the symmetric algebra $S(\hat{\mathfrak{h}}_{-})$ of $\hat{\mathfrak{h}}_{-}$
has a natural grading-restricted vertex algebra structure. In particular, 
$S(\hat{\mathfrak{h}}_{-})$ is also a 
meromorphic open-string vertex algebra. Thus 
$S(\hat{\mathfrak{h}}_{-})$
is in fact a quotient of $T(\hat{\mathfrak{h}}_{-})$ as 
meromorphic open-string vertex algebras. 
We also introduce the notion of left module for a meromorphic open-string vertex algebra
and construct a structure of left module for $T(\hat{\mathfrak{h}}_{-})$
on $T(\hat{\mathfrak{h}}_{-})\otimes M$ for a left $T(\mathfrak{h})$-module $M$. 
Comparing to the construction of the vertex operator algebra associated
to the universal enveloping algebra 
$U(\hat{\mathfrak{h}})$ of the Heisenberg algebra $\hat{\mathfrak{h}}$, the construction 
in the present paper involves much more complicated calculations because of the 
noncommutativity of the tensor algebras $T(\hat{\mathfrak{h}}_{+})$, 
$T(\hat{\mathfrak{h}}_{-})$ and $T(\mathfrak{h})$ of 
$\hat{\mathfrak{h}}_{+}$, $\hat{\mathfrak{h}}_{-}$ and $\mathfrak{h}$, respectively.

The present paper grew out of the author's study of nonlinear sigma
models using the representation theory of vertex operator algebras. 
Though there is indeed a vertex operator algebra associated to 
a Riemannian manifold and representations of this algebra can be constructed from
smooth functions on the manifold, the author has noticed that the more fundamental structure
associated to a Riemannian manifold is a meromorphic open-string vertex algebra.
It was conjectured by physicists that in general quantum nonlinear sigma models 
are not conformal field theories. These theories are believed to 
be ``gapped'' or massive theories. Therefore it is reasonable to expect that 
the fundamental algebraic structure of a nonlinear sigma model satisfies
the operator product expansion condition but in general might not have 
a commutativity property.
The construction of a meromorphic open-string vertex algebra
and left modules from a 
Riemannian manifold is given in \cite{H4}.

We also have notions of  right module and bimodule for a meromorphic 
open-string vertex algebra. We also have a similar
construction of right $T(\mathfrak{h})$-modules and we can also construct 
$T(\mathfrak{h})$-bimodules. 
These notions, constructions and a study of these modules and left modules
will be given in a paper on the representation theory of 
meromorphic 
open-string vertex algebras.

When the homogeneous subspaces of a meromorphic open-string vertex algebra is finite 
dimensional, we say that it is grading-restricted. 
The notion of grading-restricted meromorphic open-string vertex algebra should be viewed as 
a noncommutative generalization of the notion of grading-restricted vertex algebra
(which really should be called grading-restricted closed-string vertex
algebra). In \cite{H2}, the author introduced cohomologies 
of grading-restricted vertex algebras by constructing certain complexes analogous to 
the Hochschild complex for associative algebras
and then consider subcomplexes analogous to the Harrison 
complex for commutative associative algebras. The complexes in \cite{H2} that are
analogous to 
the Hochschild complex are in fact
also defined for meromorphic open-string vertex algebras
and  give cohomologies for such algebras. We shall discuss this cohomology theory in
a future publication.

The construction in the present paper can be generalized to higher dimensions.
There have been efforts by mathematicians to generalize vertex (operator) algebras
to higher dimensions. But these efforts are not very successful mainly because
the examples constructed are mostly free field theories or theories obtained 
by tensoring two-dimensional conformal field theories in a suitable sense. The main 
difficulty is that for those higher-dimensional 
quantum field theories of fundamental importance in mathematics and physics, 
there might not be commutativity 
or equivalent properties. 
On the other hand,
we do want operator product expansion to hold. Our generalizations of meromorphic 
open-string vertex algebras in higher dimensions satisfy associativity 
but not necessarily commutativity and the examples obtained by generalizing the construction 
in the present paper are not from free field theories. We shall give these 
generalizations and constructions in another future publication. 
 
The present paper is organized as follows: In Section 2, we introduce 
the notion of meromorphic open-string vertex algebra and explain that they are
indeed open-string vertex algebra defined in \cite{HK}. 
For a vector space $\mathfrak{h}$, 
we introduce a quotient algebra $N(\hat{\mathfrak{h}})$ of the tensor algebra 
$T(\hat{\mathfrak{h}})$ mentioned above 
and construct induced left modules for $N(\hat{\mathfrak{h}})$ in Section 3. 
As a preparation 
for our construction of examples of meromorphic open-string vertex algebras and left modules,
we define and study normal ordering and vertex operators in Section 4. 
This is the main technical section of the present paper. 
In Section 5, we construct the class of meromorphic open-string vertex algebras
mentioned above. In Section 6, we introduce the notion of left module for 
a meromorphic open-string vertex algebra and construct left modules for the 
meromorphic open-string vertex algebras constructed in Section 5.

\paragraph{Acknowledgments}
The author is supported in part by NSF
grant PHY-0901237.

\section{Definition of meromorphic open-string vertex algebra}

In this section, we give the definition of meromorphic open-string 
vertex algebra. We also recall the notion of open-string vertex algebra
introduced by Kong and the author in \cite{HK} and explain that 
a meromorphic open-string vertex algebra is indeed an open-string vertex algebra.

Since the applications we have in mind are always over the field of complex numbers,
for convenience, 
we shall assume that all the vector spaces in the present 
paper are over the complex numbers. But every definition, except for the 
recalled notion of open-string vertex algebra, construction or 
result in the present paper can be formulated, carried out or obtained 
over a field of characteristic $0$ without additional efforts.
We use both formal variables and complex variables. We shall use $x, x_{1}, \dots,
y, y_{1}, \dots$ to denote commuting formal variables and $z, z_{1}, \dots$ to denote
complex variables or complex numbers. When we write down the expression such as
$(x_{1}-x_{2})^{n}$ for $n\in \Z$ and commuting formal variables $x_{1}$ and $x_{2}$,
we always mean the expansion in nonnegative powers of $x_{2}$, the second formal variable.
But when we write down the expression $(z_{1}-z_{2})^{n}$ for $n\in \Z$ and complex variables
or numbers $z_{1}$ and $z_{2}$, we mean the usual analytic function or the complex number.
In the region $|z_{1}|>|z_{2}|$, this analytic function or complex number is in fact 
equal to the sum of the series obtained by substituting $z_{1}$ and $z_{2}$ 
for $x_{1}$ and $x_{2}$ in the formal series $(x_{1}-x_{2})^{n}$.

We first introduce meromorphic open-string 
vertex algebras:

\begin{defn}
{\rm A {\it meromorphic open-string vertex algebra} is a $\Z$-graded vector space 
$V=\coprod_{n\in \Z}V_{(n)}$ (graded by weights) equipped with 
a {\it vertex operator map}
\begin{eqnarray*}
Y_{V}: V&\to& ({\rm End}\; V)[[x, x^{-1}]]\nn
u&\mapsto& Y_{V}(u, x),
\end{eqnarray*}
or equivalently, 
\begin{eqnarray*}
Y_{V}: V\otimes V&\to& V[[x, x^{-1}]]\nn
u\otimes v&\mapsto& Y_{V}(u, x)v,
\end{eqnarray*}
a {\it vacuum} $\one\in V$, satisfying the 
following conditions:
\begin{enumerate}

\item {\it Lower bound condition}: When $n$ is sufficiently negative,
$V_{(n)}=0$ .

\item Properties for the vacuum: 
$Y_{V}(\one,x)=1_{V}$ (the
{\it identity property}) and for $u\in V$, $Y_{V}(u, x)\one\in V[[x]]$
and $\lim_{x\rightarrow 0} Y_{V}(u, x)\one=u$
(the {\it creation property}).

\item {\it Rationality}: For $u_{1}, \dots, u_{n}, v\in V$
and $v'\in V'$, the series 
\begin{equation}\label{product}
\langle v', Y_{V}(u_{1}, z_1)\cdots Y_{V}(u_{n}, z_n)v\rangle
\end{equation}
converges absolutely 
when $|z_1|>\cdots >|z_n|>0$ to a rational function in $z_{1}, \dots, z_{n}$
with the only possible poles at $z_{i}=0$ for $i=1, \dots, n$ and $z_{i}=z_{j}$ 
for $i\ne j$. For $u_{1}, u_{2}, v\in V$
and $v'\in V'$, the series 
\begin{equation}\label{iterate}
\langle v', Y_{V}(Y_{V}(u_{1}, z_1-z_{2})u_{2}, z_2)v\rangle
\end{equation}
converges absolutely when $|z_{2}|>|z_{1}-z_{2}|>0$ to a rational function
with the only possible poles at $z_{1}=0$, $z_{2}=0$ and $z_{1}=z_{2}$. 

\item {\it Associativity}: For $u_{1}, u_{2}, v\in V$, 
$v'\in V'$, the series 
\begin{equation}\label{associativity}
\langle v', 
Y_{V}(u_{1},z_1)Y_{V}(u_{2},z_2)v\rangle
=
\langle v', 
Y_{V}(Y_{V}(u_{1},z_{1}-z_{2})u_{2},z_2)v\rangle
\end{equation}
when  $|z_{1}|>|z_{2}|>|z_{1}-z_{2}|>0$. 

\item {\it $\mathbf{d}$-bracket property}: Let $\mathbf{d}_{V}$ be the grading 
operator on $V$, that is, $\mathbf{d}_{V}u=mu$ for $m\in \R$ and
$u\in V_{(m)}$. For $u\in V$,
\begin{equation} \label{d-com}
[\mathbf{d}_{V}, Y_{V}(u,x)]= Y_{V}(\mathbf{d}_{V}u,x)+x\frac{d}{dx}Y_{V}(u,x).
\end{equation}

\item The {\it $D$-derivative property} and the  {\it $D$-commutator formula}: 
Let $D_{V}: V\to V$ be defined by 
$$D_{V}(u)=\lim_{x\to 0}\frac{d}{dx}Y_{V}(u, x)\one$$
for $u\in V$. Then for $u\in 
V$,
\begin{eqnarray}\label{L-1}
\frac{d}{dx}Y_{V}(u, x)
&=&Y_{V}(D_{V}u, x) \nn
&=&[D_{V}, Y_{V}(u, x)].
\end{eqnarray}

\end{enumerate} 

A meromorphic open-string vertex algebra is said to be {\it grading restricted}
if $\dim V_{(n)}<\infty$ for $n\in \Z$. 
{\it Homomorphisms}, {\it isomorphisms}, {\it subalgebras}
of meromorphic open-string vertex algebras are defined in the obvious way.}
\end{defn}

We shall denote the meromorphic open-string vertex algebra defined above
by $(V, Y_{V}, \one)$ or simply by $V$. For $u\in V$, 
we call the map $Y_{V}(u, x): V\to V[[x, x^{-1}]]$ the {\it vertex 
operator associated to $u$}. 

\begin{rema}\label{op-va-nd-va}
{\rm Note that from the definition, a meromorphic open-string vertex
algebra in general might not be a vertex algebra but a $\Z$-graded 
vertex algebra such that the $\Z$-grading is lower bounded
is a meromorphic open-string vertex algebra. In particular, a grading-restricted 
vertex algebra in the sense of \cite{H2} or
a vertex operator algebra
in the sense of \cite{FLM} and \cite{FHL} is a grading-restricted meromorphic 
open-string vertex algebra. (In \cite{H2} and in the present paper,
we use the term open-string
vertex algebra
because such an algebra can be interpreted as describing the interaction of 
open stings at a "vertex." See the discussion below and the discussion in 
\cite{HK}. In fact, a (grading-restricted)
vertex algebra should have been called a (meromorphic)
closed-string vertex algebra.)}
\end{rema}

\begin{rema}
{\rm In the definition above, we require that a meromorphic open-string 
vertex algebra satisfy some strong conditions, for example, 
the lower bound condition. We can define weaker versions
of meromorphic open-string vertex algebra but here we put these stronger conditions
since the examples we construct in this paper satisfy these stronger conditions.
These conditions are also important for the development and applications of the theory of
meromorphic open-string vertex algebras.}
\end{rema}

For a $\C$-graded vector space $V=\coprod_{n\in \R}V_{(n)}$, we use $\overline{V}$
to denote the algebraic completion $\prod_{n\in \C}V_{(n)}$ of $V$.
We now recall the notion of open-string vertex algebra from \cite{HK}:

\begin{defn}
{\rm An {\it open-string vertex algebra} 
is an $\R$-graded vector space $V=\coprod_{n\in \R}V_{(n)}$ (graded
by {\it weights}) 
equipped with a {\it vertex map}
\begin{eqnarray*}
Y^{O}: V \times \R_{+} &\to& \hom(V, \overline{V})\nn
(u, r)&\mapsto &Y^{O}(u, r)
\end{eqnarray*}
or equivalently,
\begin{eqnarray*}
Y^{O}: (V\otimes V) \times \R_{+} &\to& \overline{V}\nn
(u\otimes v, r)&\mapsto &Y^{O}(u, r)v,
\end{eqnarray*}
a {\it vacuum} $\one\in V$
and an operator $D\in \edo V$ of weight $1$, satisfying 
the following conditions:

\begin{enumerate}

\item {\it Vertex map weight property}: 
For $n_{1}, n_{2}\in \R$, 
there exist a finite subset $N(n_{1}, n_{2})\subset \R$ 
such that the image of 
$\left(\coprod_{n\in n_{1}+\Z}V_{(n)}\otimes 
\coprod_{n\in n_{2}+\Z}V_{(n)}\right)\times \R_{+}$
under $Y^{O}$ is in 
$\prod_{n\in N(n_{1}, n_{2}) +\Z}
V_{(n)}.$

\item Properties for the vacuum: For any $r\in \R_{+}$,
$Y^{O}(\one,r)=1_{V}$ (the
{\it identity property}) and  $\lim_{r\rightarrow 0} Y^{O}(u, r)\one$
exists and is equal to $u$
(the {\it creation property}).

\item {\it Local-truncation property for $D'$}: Let 
$D': V'\to V'$ be the adjoint of $D$.
Then for any $v'\in V'$, there exists a positive integer $k$ such that 
$(D')^{k}v'=0$.

\item {\it Convergence properties}: For $v_{1}, \dots, v_{n}, v\in V$
and $v'\in V'$, the series 
$$\langle v', Y^{O}(v_{1}, r_1)\cdots Y^{O}(v_{n}, r_n)v\rangle$$
converges absolutely 
when $r_1>\cdots >r_n>0$. For $v_{1}, v_{2}, v\in V$
and $v'\in V'$, the series 
$$\langle v', Y^{O}(Y^{O}(v_{1}, r_0)v_{2}, r_2)v\rangle$$
converges absolutely when $r_{2}>r_{0}>0$. 

\item {\it Associativity}: For $v_{1}, v_{2}, v\in V$
and $v'\in V'$, 
$$\langle v', Y^{O}(v_{1}, r_1)Y^{O}(v_{2}, r_2)v\rangle =
\langle v', Y^{O}(Y^{O}(v_{1}, r_1-r_2)v_{2},r_2)v\rangle$$
for  $r_{1}, r_{2}\in \R$ satisfying 
$r_1>r_2>r_1-r_2>0$.

\item {\it $\mathbf{d}$-bracket property}: Let $\mathbf{d}$ be the grading 
operator on $V$, that is, $\mathbf{d}u=mu$ for $m\in \R$ and
$u\in V_{(m)}$. For $u\in V$ and $r\in \R_{+}$,
\begin{equation} \label{d-bra}
[\mathbf{d}, Y^{O}(u,r)]= Y^{O}(\mathbf{d}u,r)+r\frac{d}{dr}Y^{O}(u,r).
\end{equation}

\item {\it $D$-derivative property}: We still use $D$ to denote
the natural extension of $D$ to $\hom(\overline{V}, \overline{V})$.
For $u\in V$, 
$Y^{O}(u, r)$ as a map {}from $\R_{+}$ to 
$\hom(V, \overline{V})$
is differentiable and 
\begin{equation}\label{D-der}
\frac{d}{dr}Y^{O}(u,r)=[D, Y^{O}(u,r)]
= Y^{O}(Du,r).
\end{equation}

\end{enumerate}}
\end{defn}

The open-string vertex algebra defined above is denoted 
$(V, Y^{O}, \one, D)$ or simply $V$.

In \cite{HK}, a formal-variable vertex operator map 
\begin{eqnarray*}
\mathcal{Y}^{f}: V&\to& ({\rm End}\;V)\{x\}\nn
u&\mapsto&\Y^{f}(u, x)
\end{eqnarray*}
is constructed such that $Y^{O}(u, r)=\Y^{f}(u, r)$ for $u\in V$ and 
$r\in \R_{+}$. In particular, the open-string vertex algebra
can be studied in terms of $\Y^{f}$. 

Given a meromorphic open-string vertex algebra $(V, Y_{V}, \one)$,
let 
\begin{eqnarray*}
Y_{V}^{O}: V \times \R_{+} &\to& \hom(V, \overline{V})\nn
(u, r)&\mapsto &Y_{V}^{O}(u, r)
\end{eqnarray*}
be defined by 
$Y_{V}^{O}(u, r)=Y_{V}(u, r)$. Then we have:

\begin{prop}
The quadruple $(V, Y_{V}^{O}, \one, D_{V})$ is an open-string vertex algebra.
\end{prop}
\pf
The vertex map weight property, the identity property, the creation property,
the convergence properties, associativity, the 
$\mathbf{d}$-bracket property and the $D$-derivative property 
hold obviously. The local-truncation property for $D'$ holds because 
the meromorphic open-string vertex algebra satisfies the lower bound condition.
\epfv

In the applications of meromorphic open-string vertex algebras, we also need 
direct products of such algebras. 

\begin{defn}
{\rm Let $(V_{\alpha}, Y_{V_{\alpha}}, \mathbf{1}_{\alpha})$ 
for $\alpha\in \mathcal{A}$ be meromorphic open-string vertex algebras.
Assume that the weights of $V_{\alpha}$ for $\alpha\in \mathcal{A}$ 
are bounded from below by a common number. 
Let $V=\prod_{\alpha\in \mathcal{A}}V_{\alpha}$.
Then $V$ together with the direct products of the $\Z$-gradings, 
vertex operators and the vacuums of $V_{\alpha}$ for $\alpha\in \mathcal{A}$ 
is a meromorphic open-string vertex algebra and 
is called the {\it direct product meromorphic open-string vertex algebra}
of $(V_{\alpha}, Y_{V_{\alpha}}, \mathbf{1}_{\alpha})$,
$\alpha\in \mathcal{A}$. }
\end{defn}

\section{A quotient algebra of the tensor algebra of the affinization $\hat{\mathfrak{h}}$ 
of a vector space $\mathfrak{h}$}

Examples of open-string vertex algebras were constructed
in \cite{HK} using modules and intertwining operators for vertex operator algebras.
In this section, we study a quotient algebra
of the tensor algebra of the affinization of a vector space and its modules. We shall use
these structures in later sections to
construct directly a class of meromorphic open-string vertex operator algebras and left
modules
and thus 
new examples of open-string vertex algebras and left modules, without using the 
theory of vertex operator algebras. 

Let $\mathfrak{h}$ be a  vector space over 
$\C$ equipped with a nondegenerate bilinear form $(\cdot, \cdot)$. 
The Heisenberg algebra $\hat{\mathfrak{h}}$ associated with $\mathfrak{h}$
and $(\cdot, \cdot)$ 
is the vector space $\mathfrak{h}\otimes [t, t^{-1}]\oplus \mathbb{C}\mathbf{k}$ 
equipped with the bracket operation defined by
\begin{eqnarray*}
[a\otimes t^{m}, b\otimes t^{n}]&=&m(a,
b)\delta_{m+n, 0}\mathbf{k},\nn\\
{[a\otimes t^{m}, \mathbf{k}]}&=&0,
\end{eqnarray*}
for $a, b\in \mathfrak{h}$ and $m, n\in \mathbb{Z}$. 
It is a $\mathbb{Z}$-graded Lie algebra. 
In particular, we have the universal enveloping algebra
$U(\mathfrak{h})$ of $\mathfrak{h}$. The universal enveloping algebra
$U(\mathfrak{h})$ is constructed as a quotient of 
the tensor algebra $T(\mathfrak{h})$ of the vector space $\mathfrak{h}$.
We have a triangle decomposition
$$\hat\mathfrak{h}=\hat\mathfrak{h}_{-}\oplus \hat\mathfrak{h}_{0}
\oplus \hat\mathfrak{h}_{+},$$
where 
\begin{eqnarray*}
\hat\mathfrak{h}_{-}&=&\mathfrak{h}\otimes t^{\pm 1}\mathbb{C}[t^{-1}],\\
\hat\mathfrak{h}_{+}&=&\mathfrak{h}\otimes t^{\pm 1}\mathbb{C}[t],\\
\hat\mathfrak{h}_{0}&=&\mathfrak{h}\otimes \mathbb{C} \oplus 
\mathbb{C}\mathbf{k} \nn
&\simeq&\mathfrak{h}\oplus 
\mathbb{C}\mathbf{k},\\
\mathfrak{h}&\simeq &\mathfrak{h}\otimes \mathbb{C} 
\end{eqnarray*}
are subalgebras of 
$\hat\mathfrak{h}$.

The meromorphic open-string vertex algebras and left modules in the present paper are
constructed from left modules for a quotient algebra $N(\hat{\mathfrak{h}})$ 
of the tensor algebra $T(\hat{\mathfrak{h}})$ such that 
$U(\hat{\mathfrak{h}})$ is a quotient of $N(\hat{\mathfrak{h}})$. 
Let 
$I$ be the two-sided ideal of $T(\hat{\mathfrak{h}})$ generated by elements of the form
\begin{eqnarray*}
&(a\otimes t^{m})\otimes (b\otimes t^{n})- ((b\otimes t^{n})\otimes a\otimes t^{m})
-m(a, b)\delta_{m+n, 0}\mathbf{k},&\\
&(a\otimes t^{k})\otimes (b\otimes t^{0})-(b\otimes t^{0})\otimes (a\otimes t^{k}),&\\
&(a\otimes t^{k})\otimes \mathbf{k}-\mathbf{k}\otimes (a\otimes t^{k})&
\end{eqnarray*}
for $m\in \Z_{+}$, $n\in -\Z_{+}$, $k\in \Z$. 
Let $N(\hat{\mathfrak{h}})=T(\hat{\mathfrak{h}})/I$. 
By definition, we see that $U(\hat{\mathfrak{h}})$ is a quotient algebra
of $N(\hat{\mathfrak{h}})$. 

We have the following the Poincar\'{e}-Birkhof-Witt type
result for $U(\hat{\mathfrak{h}})$;

\begin{prop}\label{pbw-type}
As a vector space,  $N(\hat{\mathfrak{h}})$ is linearly isomorphic to 
\begin{equation}\label{pbw-type-1}
T(\hat{\mathfrak{h}}_{-})\otimes T(\hat{\mathfrak{h}}_{+})\otimes T(\mathfrak{h})
\otimes T(\mathbb{C}\mathbf{k})
\end{equation}
where $T(\hat{\mathfrak{h}}_{-})$, $T(\hat{\mathfrak{h}}_{+})$, $T(\mathfrak{h})$
and $T(\mathbb{C}\mathbf{k})$
are the tensor algebras of the vector spaces $\hat{\mathfrak{h}}_{-}$,
$\hat{\mathfrak{h}}_{+}$, $\mathfrak{h}$ and $\mathbb{C}\mathbf{k}$, respectively.
\end{prop}
\pf
We first show that for any $k\in \N$,
any element of $T(\hat{\mathfrak{h}})$ of the form $u_{1}\otimes
\cdots \otimes u_{k}$ for $u_{1}, \dots, u_{k}$ of either the form $a\otimes t^{m}$ or 
$\mathbf{k}$ is a sum of an element of (\ref{pbw-type-1})
and an element of $I$. We use induction on the number of elements that are not 
$\mathbf{k}$ in the set $\{u_{1}, \dots, u_{k}\}$. When this number is $0$, 
the element we are considering is 
$\mathbf{k}\otimes
\cdots \otimes \mathbf{k}$ and is in
(\ref{pbw-type-1}).
Assume that when there are less than $n$ elements that are not 
$\mathbf{k}$ in the set $\{u_{1}, \dots, u_{k}\}$, this statement is true. 
Modulo elements of $I$, 
we can move all factors of the form $\mathbf{k}$ in the element $u_{1}\otimes
\cdots \otimes u_{k}$
to the right and then move those $u_{i}\in \mathfrak{h}_{0}$ to the immediate left of 
the tensor powers of $\mathbf{k}$ but keep the order of these elements. 
Thus we can assume that $u_{1}\otimes
\cdots \otimes u_{k}$ is of the form 
$$a_{1}\otimes t^{m_{1}}\otimes \cdots a_{l}\otimes t^{m_{l}}
\otimes a_{l+1}\otimes \cdots\otimes a_{n}\otimes \mathbf{k}\otimes 
\cdots\otimes \mathbf{k}$$
where $a_{1}, \dots, a_{n}\in \mathfrak{h}$ and $m_{1}, 
\dots, m_{l}\in \Z\setminus \{0\}$. If there is an integer $j$ satisfying 
$1\le j\le l$ such that $m_{1}, \dots, m_{j}<0$ and $m_{j+1}, \dots, m_{l}>0$,
then this element is in (\ref{pbw-type-1}).
Otherwise, modulo elements of $I$ and elements of the form $u_{1}\otimes
\cdots \otimes u_{k}$ with less than $n$ factors 
not equal to $\mathbf{k}$, we can move factors of the form 
$a_{i}\otimes t^{m_{i}}$ with positive $m_{i}$ to the right of the factors of 
the form $a_{i}\otimes t^{m_{i}}$ with negative $m_{i}$ and keep the order
of such factors with positive $m_{i}$ and the order of such factors 
with negative $m_{i}$. The resulting element is in (\ref{pbw-type-1}).
By induction assumption, elements of the form $u_{1}\otimes
\cdots \otimes u_{k}$ with less than $n$ factors 
not equal to $\mathbf{k}$ are sums of elements of (\ref{pbw-type-1})
and $I$. Thus in the case that there are $n$ elements that are not 
$\mathbf{k}$ in the set $\{u_{1}, \dots, u_{k}\}$, the statement is true.

We have proved that $T(\hat{\mathfrak{h}})$ is the sum of (\ref{pbw-type-1})
and $I$. Since the intersection of  (\ref{pbw-type-1}) and $I$ is clearly $0$, 
$T(\hat{\mathfrak{h}})$ is the direct sum of (\ref{pbw-type-1})
and $I$. Thus $N(\hat{\mathfrak{h}})=T(\hat{\mathfrak{h}})/I$ is linearly isomorphic to 
(\ref{pbw-type-1}).
\epfv

Now we construct left modules for $N(\hat{\mathfrak{h}})$. Let 
$M$ be a left $T(\mathfrak{h})$-module. We define the action of 
$\mathbf{k}$ on $M$ to be $1$ and the actions of elements of $\mathfrak{h}_{+}$
on $M$ to be $0$. Then $M$ is also
a left module for the subalgebra $N(\hat{\mathfrak{h}}_{+}\oplus \hat{\mathfrak{h}}_{0})$ 
of $N(\hat{\mathfrak{h}})$ generated by
elements of $\hat{\mathfrak{h}}_{+}$ and $\hat{\mathfrak{h}}_{0}$. 
We consider the induced left module 
$N(\hat{\mathfrak{h}})\otimes_{N(\hat{\mathfrak{h}}_{+}\oplus \hat{\mathfrak{h}}_{0})} M$.
By Proposition \ref{pbw-type}, we see that 
$N(\hat{\mathfrak{h}})\otimes_{N(\hat{\mathfrak{h}}_{+}\oplus \hat{\mathfrak{h}}_{0})} M$ 
is linearly isomorphic to $T(\hat{\mathfrak{h}}_{-})\otimes M$. We shall identify 
$N(\hat{\mathfrak{h}})\otimes_{N(\hat{\mathfrak{h}}_{+}\oplus \hat{\mathfrak{h}}_{0})} M$ 
with $T(\hat{\mathfrak{h}}_{-})\otimes M$. The left $N(\hat{\mathfrak{h}})$-module structure
on $T(\hat{\mathfrak{h}}_{-})\otimes M$ can be obtained explicitly by using the commutator
relations defining the algebra $N(\hat{\mathfrak{h}})$ and the left 
$N(\hat{\mathfrak{h}}_{+}\oplus  \hat{\mathfrak{h}}_{0})$-module structure on $M$. 

For a left $N(\hat{\mathfrak{h}})$-module,
we denote the representation images of $a\otimes t^{n}\in \hat{\mathfrak{h}}$ 
for $a\in \mathfrak{h}$ and $n\in \Z$ acting on the left module by $a(n)$. Then 
a left $N(\hat{\mathfrak{h}})$-module
$T(\hat{\mathfrak{h}}_{-})\otimes M$ constructed from a left $T(\mathfrak{h})$-module $M$ 
is spanned by elements of the form
$a_{1}(-n_{1})\cdots a_{k}(-n_{k})w$, where 
$a_{1}, \dots, a_{k}\in \hat{\mathfrak{h}}$, $n_{1}, \dots, n_{k}\in \Z_{+}$
and $w\in M$.

\section{Normal ordering and vertex operators}

In this section, we define the normal ordering for certain operators on 
a left $N(\hat{\mathfrak{h}})$-module of the form $T(\hat{\mathfrak{h}}_{-})\otimes M$
and vertex operators 
acting on such a left $N(\hat{\mathfrak{h}})$-module.
We then prove a number of technical formulas for products of normal ordered 
products of operators and products of vertex operators. This section contain the 
main technical material of the present paper. Many of the calculations are much more
complicated than the Heisenberg algebra case 
because of the noncommutativity of the operators.

Given a left $N(\hat{\mathfrak{h}})$-module,
we define a {\it normal ordering} map $\nord \cdot \nord$
from the space of operators on 
the left module spanned by operators of the form $a_{1}(n_{1})\cdots a_{k}(n_{k})$
to itself by 
$$\nord a_{1}(n_{1})\cdots a_{k}(n_{k})\nord=a_{\sigma(1)}(n_{\sigma(1)})
\cdots a_{\sigma(k)}(n_{\sigma(k)}),$$
where $\sigma\in S_{k}$ is the unique permutation such that 
\begin{eqnarray*}
&\sigma(1)<\cdots <\sigma(\alpha),&\\ 
&\sigma(\alpha+1)<\cdots <\sigma(\beta),&\\
&\sigma(\beta+1)<\cdots <\sigma(k),&\\
&n_{\sigma(1)}, \dots, n_{\sigma(\alpha)}<0,&\\
&n_{\sigma(\alpha+1)}, \dots, n_{\sigma(\beta)}>0,&\\
&n_{\sigma(\beta+1)}, \dots, n_{\sigma(k)}=0,&
\end{eqnarray*}
for some integers $\alpha$ and $\beta$ satisfying $0\le \alpha\le \beta\le k$.

Given an induced left $N(\hat{\mathfrak{h}})$-module $W=T(\hat{\mathfrak{h}}_{-})\otimes M$,
$a_{1}, \dots, a_{k}\in
\mathfrak{h}$ and $m_{1}, \dots, m_{k}\in \mathbb{Z}_{+}$, we define
the vertex operator $Y_{W}(a_{1}(-m_{1})\cdots a_{k}(-m_{k})\one, x)$
associated to $a_{1}(-m_{1})\cdots a_{k}(-m_{k})\one\in T(\hat{\mathfrak{h}}_{-})$
by
\begin{eqnarray}\label{vo-map}
\lefteqn{Y_{W}(a_{1}(-m_{1})\cdots a_{k}(-m_{k})\one, x)}\nn
&&=\nord \frac{1}{(n_{1}-1)!}\left(\frac{d^{m_{1}-1}}{dx^{m_{1}-1}}
a_{1}(x)\right)\cdots \frac{1}{(m_{k}-1)!}
\left(\frac{d^{m_{k}-1}}{dx^{m_{k}-1}}
a_{k}(x)\right)\nord,
\end{eqnarray}
where 
$$a_{i}(x)=\sum_{n\in \Z}a_{i}(n)x^{-n-1}$$
for $i=1, \dots, k$ and $a_{i}(n)$ for $i=1, \dots, k$ and $n\in \Z$ are the 
representation images of $a_{i}\otimes t^{n}$ on $W$. 

We need the following commutator formula:

\begin{lemma}
For $a, b\in \mathfrak{h}$, 
\begin{eqnarray}\label{commutator-pm}
\lefteqn{\left[\frac{1}{(m-1)!}\frac{\partial^{m-1}}{\partial x_{1}^{m-1}}
a^{+}(x_{1}), \frac{1}{(n-1)!}\frac{\partial^{n-1}}{\partial x_{2}^{n-1}}
b^{-}(x_{2})\right]}\nn
&&=n(a, b){-n-1\choose m-1}(x_{1}-x_{2})^{-m-n},
\end{eqnarray}
where for $a\in \mathfrak{h}$, 
$$a^{\pm}(x)=\sum_{n\in \pm \Z_{+}}a(n)x^{-n-1}$$
and a negative power of $x_{1}-x_{2}$, as in the formal calculus in the theory of vertex 
operator algebras, is understood as the binomial expansion in the nonnegative 
powers of the formal variable $x_{2}$.
\end{lemma}
\pf
The proof is a straightforward calculation.
\epfv

We also need an explicit expression of a vertex operator. 
For $k\in \Z_{+}$ and $\alpha, \beta\in \N$ satisfying
$0\le \alpha\le \beta\le k$, let $J(k; \alpha, \beta)$
be the set of elements of 
$S_{k}$ which preserve the orders of the first $\alpha$ numbers, the next
$\beta-\alpha$ numbers, and the last $k-\beta$ numbers, that is,
\begin{eqnarray*}
\lefteqn{J_{k; \alpha, \beta}=\{\sigma\in S_{k}\;|\;\sigma(1)<\cdots <\sigma(\alpha),\;}\nn
&&\quad\quad\quad\quad\quad\quad\quad
\sigma(\alpha+1)<\cdots <\sigma(\beta),\; \sigma(\beta+1)<\cdots <\sigma(k)\}.
\end{eqnarray*}

\begin{lemma}
For $a_{1}, \dots, a_{k}\in \hat{\mathfrak{h}}$ and $n_{1}, \dots, n_{k}\in \Z_{+}$, 
\begin{eqnarray}\label{vo-exp}
\lefteqn{Y_{W}(a_{1}(-n_{1})\cdots a_{k}(-n_{k})\one, x)}\nn
&&=\sum_{0\le \alpha\le \beta\le k}\sum_{\sigma\in J(k; \alpha, \beta)}
\left(\frac{1}{(n_{\sigma(1)}-1)!}\frac{\partial^{n_{\sigma(1)}-1}}
{\partial z_{2}^{n_{\sigma(1)}-1}}
a^{-}_{\sigma(1)}(z_{2})\right)\cdot\nn
&&\quad\quad\quad\quad\quad\quad\quad\quad\quad \cdots
\left(\frac{1}{(n_{\sigma(\alpha)}-1)!}\frac{\partial^{n_{\sigma(\alpha)}-1}}
{\partial z_{2}^{n_{\sigma(\alpha)}-1}}
a^{-}_{\sigma(\alpha)}(z_{2})\right)\cdot\nn
&&\quad\quad\quad\quad\quad\quad\quad\quad\cdot \left(\frac{1}{(n_{\sigma(\alpha+1)}-1)!}
\frac{\partial^{n_{\sigma(\alpha+1)}-1}}{\partial z_{2}^{n_{\sigma(\alpha+1)}-1}}
a^{+}_{\sigma(\alpha+1)}(z_{2})\right)\cdot\nn
&&\quad\quad\quad\quad\quad\quad\quad\quad\quad\cdots
\left(\frac{1}{(n_{\sigma(\beta)}-1)!}\frac{\partial^{n_{\sigma(\beta)}-1}}
{\partial z_{2}^{n_{\sigma(\beta)}-1}}
a^{+}_{\sigma(\beta)}(z_{2})\right)\cdot\nn
&&\quad\quad\quad\quad\quad\quad\quad\quad\cdot \left(\frac{1}{(n_{\sigma(\alpha+1)}-1)!}
\frac{\partial^{n_{\sigma(\beta+1)}-1}}{\partial z_{2}^{n_{\sigma(\beta+1)}-1}}
a_{\sigma(\beta+1)}(0)z_{2}^{-1}\right)\cdot\nn
&&\quad\quad\quad\quad\quad\quad\quad\quad\quad\cdots
\left(\frac{1}{(n_{\sigma(k)}-1)!}\frac{\partial^{n_{\sigma(k)}-1}}
{\partial z_{2}^{n_{\sigma(k)}-1}}
(a_{\sigma(k)}(0)z_{2}^{-1})\right).
\end{eqnarray}
\end{lemma}
\pf
The expression follows immediately from the definition of the 
vertex operator and the definition of the normal ordering.
\epfv

We need the following:

\begin{prop}
For $a_{0}, \dots, a_{k}\in \mathfrak{h}$ and $m_{0}, \dots, m_{k}\in \Z_{+}$, 
\begin{eqnarray}\label{normal-ordering}
\lefteqn{\left(\frac{1}{(m_{0}-1)!}\frac{\partial^{m_{0}-1}}
{\partial x_{0}^{m_{0}-1}}a_{0}(x_{0})\right)\cdot}\nn
&&\quad\quad\cdot
\nord \left(\frac{1}{(m_{1}-1)!}\frac{\partial^{m_{1}-1}}
{\partial x_{1}^{m_{1}-1}}a_{1}(x_{1})\right)
\cdots \left(\frac{1}{(m_{k}-1)!}\frac{\partial^{m_{k}-1}}
{\partial x_{k}^{m_{k}-1}}a_{k}(x_{k})\right)\nord\nn
&&=\nord \left(\frac{1}{(m_{0}-1)!}\frac{\partial^{m_{0}-1}}
{\partial x_{0}^{m_{0}-1}}a_{0}(x_{0})\right)\cdot\nn
&&\quad\quad\cdot
\left(\frac{1}{(m_{1}-1)!}\frac{\partial^{m_{1}-1}}
{\partial x_{1}^{m_{1}-1}}a_{1}(x_{1})\right)
\cdots \left(\frac{1}{(m_{k}-1)!}\frac{\partial^{m_{k}-1}}
{\partial x_{k}^{m_{k}-1}}a_{k}(x_{k})\right)\nord\nn
&&\quad +\sum_{p=1}^{k}m_{p}(a_{0}, a_{p}){-m_{p}-1\choose m_{0}-1}
(x_{0}-x_{p})^{-m_{0}-m_{p}}\cdot\nn
&&\quad\quad\cdot \nord \left(\frac{1}{(m_{1}-1)!}\frac{\partial^{m_{1}-1}}
{\partial x_{1}^{m_{1}-1}}a_{1}(x_{1})\right)\cdot\nn
&&\quad\quad\quad\quad\cdots
\left(\frac{1}{(m_{p-1}-1)!}\frac{\partial^{m_{p-1}-1}}
{\partial x_{p-1}^{m_{p-1}-1}}a_{p-1}(x_{p-1})\right)\cdot\nn
&&\quad\quad\quad\quad\cdot \left(\frac{1}{(m_{p+1}-1)!}\frac{\partial^{m_{p+1}-1}}
{\partial x_{p+1}^{m_{p+1}-1}}a_{p+1}(x_{p+1})\right)\cdot\nn
&&\quad\quad\quad\quad
\cdots \left(\frac{1}{(m_{k}-1)!}\frac{\partial^{m_{k}-1}}
{\partial x_{k}^{m_{k}-1}}a_{k}(x_{k})\right)\nord.
\end{eqnarray}
\end{prop}
\pf
The proof is a tedious but straightforward calculation. 
By (\ref{vo-exp}) and (\ref{commutator-pm}), 
the left-hand side of (\ref{normal-ordering}) is equal to 
\begin{eqnarray*}
\lefteqn{\left(\frac{1}{(m_{0}-1)!}\frac{\partial^{m_{0}-1}}{\partial x_{0}^{m_{0}-1}}
(a_{0}^{+}(x_{0})+a_{0}(0)x_{0}^{-1}+a_{0}^{-}(x_{0}))\right)\cdot}\nn
&&\quad\quad\quad\cdot 
\left(\sum_{0\le \alpha\le \beta\le k}\sum_{\sigma\in J(k; \alpha, \beta)}
\left(\frac{1}{(m_{\sigma(1)}-1)!}\frac{\partial^{m_{\sigma(1)}-1}}{\partial 
x_{\sigma(1)}^{m_{\sigma(1)}-1}}
a^{-}_{\sigma(1)}(x_{\sigma(1)})\right)\cdot\right.\nn
&&\quad\quad\quad\quad\quad\quad \cdots
\left(\frac{1}{(m_{\sigma(\alpha)}-1)!}\frac{\partial^{m_{\sigma(\alpha)}-1}}
{\partial x_{\sigma(\alpha)}^{m_{\sigma(\alpha)}-1}}
a^{-}_{\sigma(\alpha)}(x_{\sigma(\alpha)})\right)\cdot\nn
&&\quad\quad\quad\quad\quad\cdot \left(\frac{1}{(m_{\sigma(\alpha+1)}-1)!}
\frac{\partial^{m_{\sigma(\alpha+1)}-1}}{\partial x_{\sigma(\alpha+1)}^{m_{\sigma(\alpha+1)}-1}}
a^{+}_{\sigma(\alpha+1)}(x_{\sigma(\alpha+1)})\right)\cdot\nn
&&\quad\quad\quad\quad\quad\quad\cdots
\left(\frac{1}{(m_{\sigma(\beta)}-1)!}\frac{\partial^{m_{\sigma(\beta)}-1}}
{\partial x_{\sigma(\beta)}^{m_{\sigma(\beta)}-1}}
a^{+}_{\sigma(\beta)}(x_{\sigma(\beta)})\right)\cdot\nn
&&\quad\quad\quad\quad\quad\cdot \left(\frac{1}{(m_{\sigma(\alpha+1)}-1)!}
\frac{\partial^{m_{\sigma(\beta+1)}-1}}{\partial x_{\sigma(\alpha+1)}^{m_{\sigma(\beta+1)}-1}}
a_{\sigma(\beta+1)}(0)x_{\sigma(\alpha+1)}^{-1}\right)\cdot\nn
&&\quad\quad\quad\quad\quad\quad\cdots
\left.\left(\frac{1}{(m_{\sigma(l)}-1)!}\frac{\partial^{m_{\sigma(k)}-1}}
{\partial x_{\sigma(k)}^{m_{\sigma(k)}-1}}
(a_{\sigma(k)}(0)x_{\sigma(k)}^{-1})\right)\right)\nn
&&=\nord\left(\frac{1}{(m_{0}-1)!}\frac{\partial^{m_{0}-1}}{\partial x_{0}^{m_{0}-1}}
a_{0}(x_{0})\right)\cdot\nn
&&\quad\quad\quad\cdot  \left(\frac{1}{(m_{1}-1)!}\frac{\partial^{m_{1}-1}}
{\partial x_{1}^{m_{1}-1}}a_{1}(x_{1})\right)
\cdots \left(\frac{1}{(m_{k}-1)!}\frac{\partial^{m_{k}-1}}
{\partial x_{k}^{m_{k}-1}}a_{k}(x_{k})\right)\nord \nn
&&\quad +\sum_{0\le \alpha\le \beta\le k}
\sum_{\sigma\in J(k; \alpha, \beta)}\sum_{p=1}^{k}\nn
&&\quad\quad\left(
\left(\frac{1}{(m_{\sigma(1)}-1)!}\frac{\partial^{m_{\sigma(1)}-1}}
{\partial x_{\sigma(1)}^{m_{\sigma(1)}-1}}
a^{-}_{\sigma(1)}(x_{\sigma(1)})\right)\cdot\right.\nn
&&\quad\quad\quad\quad\quad\quad \cdots
\left(\frac{1}{(m_{\sigma(p-1)}-1)!}\frac{\partial^{m_{\sigma(p-1)}-1}}
{\partial x_{\sigma(p-1)}^{m_{\sigma(p-1)}-1}}
a^{-}_{\sigma(p-1)}(x_{\sigma(p-1)})\right)\cdot\nn
&&\quad\quad\cdot
\left[\left(\frac{1}{(m_{0}-1)!}\frac{\partial^{m_{0}-1}}{\partial x_{0}^{m_{0}-1}}
a_{0}^{+}(x_{0})\right), \left(\frac{1}{(m_{\sigma(p)}-1)!}\frac{\partial^{m_{\sigma(p)}-1}}
{\partial x_{\sigma(p)}^{m_{\sigma(p)}-1}}
a^{-}_{\sigma(p)}(x_{\sigma(p)})\right)\right]\cdot\nn
&&\quad \quad\quad \quad\quad\cdot
\left(\frac{1}{(m_{\sigma(p+1)}-1)!}\frac{\partial^{m_{\sigma(p+1)}-1}}
{\partial x_{\sigma(p+1)}^{m_{\sigma(p+1)}-1}}
a^{-}_{\sigma(p+1)}(x_{\sigma(p+1)})\right)\cdot\nn
&&\quad\quad\quad\quad\quad\quad \cdots
\left(\frac{1}{(m_{\sigma(\alpha)}-1)!}\frac{\partial^{m_{\sigma(\alpha)}-1}}
{\partial x_{\sigma(\alpha)}^{m_{\sigma(\alpha)}-1}}
a^{-}_{\sigma(\alpha)}(x_{\sigma(\alpha)})\right)\cdot\nn
&&\quad\quad\quad\quad\quad\cdot \left(\frac{1}{(m_{\sigma(\alpha+1)}-1)!}
\frac{\partial^{m_{\sigma(\alpha+1)}-1}}{\partial x_{\sigma(\alpha+1)}^{m_{\sigma(\alpha+1)}-1}}
a^{+}_{\sigma(\alpha+1)}(x_{\sigma(\alpha+1)})\right)\cdot\nn
&&\quad\quad\quad\quad\quad\quad\cdots
\left(\frac{1}{(m_{\sigma(\beta)}-1)!}\frac{\partial^{m_{\sigma(\beta)}-1}}
{\partial x_{\sigma(\beta)}^{m_{\sigma(\beta)}-1}}
a^{+}_{\sigma(\beta)}(x_{\sigma(\beta)})\right)\cdot\nn
&&\quad\quad\quad\quad\quad\cdot \left(\frac{1}{(m_{\sigma(\alpha+1)}-1)!}
\frac{\partial^{m_{\sigma(\beta+1)}-1}}{\partial x_{\sigma(\alpha+1)}^{m_{\sigma(\beta+1)}-1}}
a_{\sigma(\beta+1)}(0)x_{\sigma(\alpha+1)}^{-1}\right)\cdot\nn
&&\quad\quad\quad\quad\quad\quad\cdots
\left.\left(\frac{1}{(m_{\sigma(k)}-1)!}\frac{\partial^{m_{\sigma(k)}-1}}
{\partial x_{\sigma(k)}^{m_{\sigma(k)}-1}}
a_{\sigma(k)}(0)x_{\sigma(k)}^{-1})\right)\right) \nn
&&=\nord\left(\frac{1}{(m_{0}-1)!}\frac{\partial^{m_{0}-1}}{\partial x_{0}^{m_{0}-1}}
a_{0}(x_{0})\right)\cdot\nn
&&\quad\quad\quad\cdot  \left(\frac{1}{(m_{1}-1)!}\frac{\partial^{m_{1}-1}}
{\partial x_{1}^{m_{1}-1}}a_{1}(x_{1})\right)
\cdots \left(\frac{1}{(m_{k}-1)!}\frac{\partial^{m_{k}-1}}
{\partial x_{k}^{m_{k}-1}}a_{k}(x_{k})\right)\nord\nn
&&\quad +\sum_{0\le \alpha\le \beta\le k}\sum_{\sigma\in J(k; \alpha, \beta)}\sum_{p=1}^{k}
m_{\sigma(p)}(a_{0}, a_{\sigma(p)}){-m_{\sigma(p)}-1\choose m_{0}-1}
(x_{0}-x_{\sigma(p)})^{-m_{0}-m_{\sigma(p)}}\cdot\nn
&&\quad \quad \cdot\left(
\left(\frac{1}{(m_{\sigma(1)}-1)!}\frac{\partial^{m_{\sigma(1)}-1}}
{\partial x_{\sigma(1)}^{m_{\sigma(1)}-1}}
a^{-}_{\sigma(1)}(x_{\sigma(1)})\right)\cdot\right.\nn
&&\quad\quad\quad\quad\quad\quad \cdots
\left(\frac{1}{(m_{\sigma(p-1)}-1)!}\frac{\partial^{m_{\sigma(p-1)}-1}}
{\partial x_{\sigma(p-1)}^{m_{\sigma(p-1)}-1}}
a^{-}_{\sigma(p-1)}(x_{\sigma(p-1)})\right)\cdot\nn
&&\quad \quad\quad \quad\quad\cdot
\left(\frac{1}{(m_{\sigma(p+1)}-1)!}\frac{\partial^{m_{\sigma(p+1)}-1}}
{\partial x_{\sigma(p+1)}^{m_{\sigma(p+1)}-1}}
a^{-}_{\sigma(p+1)}(x_{\sigma(p+1)})\right)\cdot\nn
&&\quad\quad\quad\quad\quad\quad \cdots
\left(\frac{1}{(m_{\sigma(\alpha)}-1)!}\frac{\partial^{m_{\sigma(\alpha)}-1}}
{\partial x_{\sigma(\alpha)}^{m_{\sigma(\alpha)}-1}}
a^{-}_{\sigma(\alpha)}(x_{\sigma(\alpha)})\right)\cdot\nn
&&\quad\quad\quad\quad\quad\cdot \left(\frac{1}{(m_{\sigma(\alpha+1)}-1)!}
\frac{\partial^{m_{\sigma(\alpha+1)}-1}}{\partial x_{\sigma(\alpha+1)}^{m_{\sigma(\alpha+1)}-1}}
a^{+}_{\sigma(\alpha+1)}(x_{\sigma(\alpha+1)})\right)\cdot\nn
&&\quad\quad\quad\quad\quad\quad\cdots
\left(\frac{1}{(m_{\sigma(\beta)}-1)!}\frac{\partial^{m_{\sigma(\beta)}-1}}
{\partial x_{\sigma(\beta)}^{m_{\sigma(\beta)}-1}}
a^{+}_{\sigma(\beta)}(x_{\sigma(\beta)})\right)\cdot\nn
&&\quad\quad\quad\quad\quad\cdot \left(\frac{1}{(m_{\sigma(\alpha+1)}-1)!}
\frac{\partial^{m_{\sigma(\beta+1)}-1}}{\partial x_{\sigma(\beta+1)}^{m_{\sigma(\beta+1)}-1}}
a_{\sigma(\beta+1)}(0)x_{\sigma(\beta+1)}^{-1}\right)\cdot\nn
&&\quad\quad\quad\quad\quad\quad\cdots
\left.\left(\frac{1}{(m_{\sigma(k)}-1)!}\frac{\partial^{m_{\sigma(k)}-1}}
{\partial x_{\sigma(k)}^{m_{\sigma(k)}-1}}
(a_{\sigma(k)}(0)x_{\sigma(k)}^{-1})\right)\right)\nn
&&=\nord\left(\frac{1}{(m_{0}-1)!}\frac{\partial^{m_{0}-1}}{\partial x_{0}^{m_{0}-1}}
a_{0}(x_{0})\right)\cdot\nn
&&\quad\quad\quad\cdot  \left(\frac{1}{(m_{1}-1)!}\frac{\partial^{m_{1}-1}}
{\partial x_{1}^{m_{1}-1}}a_{1}(x_{1})\right)
\cdots \left(\frac{1}{(m_{k}-1)!}\frac{\partial^{m_{k}-1}}
{\partial x_{k}^{m_{k}-1}}a_{k}(x_{k})\right)\nord \nn
&&\quad +\sum_{0\le \alpha\le \beta\le k}\sum_{\sigma\in J(k; \alpha, \beta)}\sum_{p=1}^{k}
m_{p}(a_{0}, a_{p}){-m_{p}-1\choose m_{0}-1}(x_{0}-x_{p})^{-m_{0}-m_{p}}\cdot\nn
&&\quad \quad\cdot\left(
\left(\frac{1}{(m_{\sigma(1)}-1)!}\frac{\partial^{m_{\sigma(1)}-1}}
{\partial x_{1}^{m_{\sigma(1)}-1}}
a^{-}_{\sigma(1)}(x_{1})\right)\cdot\right.\nn
&&\quad\quad\quad\quad\quad\quad \cdots
\left(\frac{1}{(m_{\sigma(p-1)}-1)!}\frac{\partial^{m_{\sigma(p-1)}-1}}
{\partial x_{\sigma(p-1)}^{m_{\sigma(p-1)}-1}}
a^{-}_{\sigma(p-1)}(x_{\sigma(p-1)})\right)\cdot\nn
&&\quad \quad\quad \quad\quad\cdot
\left(\frac{1}{(m_{\sigma(p+1)}-1)!}\frac{\partial^{m_{\sigma(p+1)}-1}}
{\partial x_{\sigma(p+1)}^{m_{\sigma(p+1)}-1}}
a^{-}_{\sigma(p+1)}(x_{\sigma(p+1)})\right)\cdot\nn
&&\quad\quad\quad\quad\quad\quad \cdots
\left(\frac{1}{(m_{\sigma(\alpha)}-1)!}\frac{\partial^{m_{\sigma(\alpha)}-1}}
{\partial x_{\sigma(\alpha)}^{m_{\sigma(\alpha)}-1}}
a^{-}_{\sigma(\alpha)}(x_{\sigma(\alpha)})\right)\cdot\nn
&&\quad\quad\quad\quad\quad\cdot \left(\frac{1}{(m_{\sigma(\alpha+1)}-1)!}
\frac{\partial^{m_{\sigma(\alpha+1)}-1}}{\partial x_{\sigma(\alpha+1)}^{m_{\sigma(\alpha+1)}-1}}
a^{+}_{\sigma(\alpha+1)}(x_{\sigma(\alpha+1)})\right)\cdot\nn
&&\quad\quad\quad\quad\quad\quad\cdots
\left(\frac{1}{(m_{\sigma(\beta)}-1)!}\frac{\partial^{m_{\sigma(\beta)}-1}}
{\partial x_{\sigma(\beta)}^{m_{\sigma(\beta)}-1}}
a^{+}_{\sigma(\beta)}(x_{\sigma(\beta)})\right)\cdot\nn
&&\quad\quad\quad\quad\quad\cdot \left(\frac{1}{(m_{\sigma(\alpha+1)}-1)!}
\frac{\partial^{m_{\sigma(\beta+1)}-1}}{\partial x_{\sigma(\beta+1)}^{m_{\sigma(\beta+1)}-1}}
a_{\sigma(\beta+1)}(0)x_{\sigma(\beta+1)}^{-1}\right)\cdot\nn
&&\quad\quad\quad\quad\quad\quad\cdots
\left.\left(\frac{1}{(m_{\sigma(k)}-1)!}\frac{\partial^{m_{\sigma(k)}-1}}
{\partial x_{\sigma(k)}^{m_{\sigma(k)}-1}}
(a_{\sigma(k)}(0)x_{\sigma(k)}^{-1})\right)\right) \nn
&&=\nord\left(\frac{1}{(m_{0}-1)!}\frac{\partial^{m_{0}-1}}{\partial x_{0}^{m_{0}-1}}
a_{0}(x_{0})\right)\cdot\nn
&&\quad\quad\quad\cdot  \left(\frac{1}{(m_{1}-1)!}\frac{\partial^{m_{1}-1}}
{\partial x_{1}^{m_{1}-1}}a_{1}(x_{1})\right)
\cdots \left(\frac{1}{(m_{k}-1)!}\frac{\partial^{m_{k}-1}}
{\partial x_{k}^{m_{k}-1}}a_{k}(x_{k})\right)\nord \nn
&&\quad +\sum_{p=1}^{k}m_{p}(a_{0}, a_{p}){-m_{p}-1\choose m_{0}-1}(x_{0}-x_{p})^{-m_{0}-m_{p}}
\cdot\nn
&&\quad \quad \cdot \nord\left(\frac{1}{(m_{1}-1)!}\frac{\partial^{m_{1}-1}}
{\partial x_{1}^{m_{1}-1}}a_{1}(x_{1})\right)\cdot\nn
&&\quad\quad\quad\quad\quad\cdots
\left(\frac{1}{(m_{p-1}-1)!}\frac{\partial^{m_{p-1}-1}}
{\partial x_{p-1}^{m_{p-1}-1}}a_{p-1}(x_{p-1})\right)\cdot\nn
&&\quad \quad \quad\quad\quad\quad\cdot 
\left(\frac{1}{(m_{p+1}-1)!}\frac{\partial^{m_{p+1}-1}}
{\partial x_{p+1}^{m_{p+1}-1}}a_{p+1}(x_{p+1})\right)\cdot\nn
&&\quad \quad \quad\quad\quad\quad\quad\quad
\cdots \left(\frac{1}{(m_{k}-1)!}\frac{\partial^{m_{k}-1}}
{\partial x_{k}^{m_{k}-1}}a_{k}(x_{k})\right)\nord.
\end{eqnarray*}
\epfv

\begin{rema}\label{formal-limit}
{\rm In (\ref{normal-ordering}),  the formal 
variables $x_{1}, \dots, x_{k}$ can be taken to be the same but cannot
be equal to $x_{0}$.}
\end{rema}

From (\ref{normal-ordering}), we obtain immediately the following result:

\begin{cor}
For $a_{0}, \dots, a_{k}\in \mathfrak{h}$ and $m_{0}, \dots, m_{k}\in \Z_{+}$, 
\begin{eqnarray}\label{k=1}
\lefteqn{Y_{W}(a_{0}(-m_{0})\one,x_1)
Y_{W}(a_{1}(-m_{1})\cdots a_{k}(-m_{k})\one, x_2)}\nn
&&=\nord\left(\frac{1}{(m_{0}-1)!}\frac{\partial^{m_{0}-1}}{\partial x_{1}^{m_{0}-1}}
a_{0}(x_{1})\right)\cdot\nn
&&\quad\quad\quad\cdot  \left(\frac{1}{(m_{1}-1)!}\frac{\partial^{m_{1}-1}}
{\partial x_{2}^{m_{1}-1}}a_{1}(x_{2})\right)
\cdots \left(\frac{1}{(m_{k}-1)!}\frac{\partial^{m_{k}-1}}
{\partial x_{2}^{m_{k}-1}}a_{k}(x_{2})\right)\nord \nn
&&\quad +\sum_{p=1}^{l}m_{p}(a_{0}, a_{p}){-m_{p}-1\choose m_{0}-1}
(x_{1}-x_{2})^{-m_{0}-m_{p}}\cdot
\nn
&&\quad \quad \quad\cdot Y_{W}(a_{1}(-m_{1})\cdots \widehat{a_{1}(-m_{p})}
\cdots a_{k}(-m_{k})\one,x_2),
\end{eqnarray}
where for $p=1, \dots, k$, 
we use $\widehat{a_{p}(-m_{p})}$ to denote that $a_{p}(-m_{p})$ is missing 
from a product.
\end{cor}
\pf
Taking $x_{0}$ to be $x_{1}$ and $x_{1}, \dots, x_{k} $ to be $x_{2}$
in (\ref{normal-ordering}) (note Remark \ref{formal-limit}), we obtain (\ref{k=1}). 
\epfv

\begin{cor}
For $a_{0}, \dots, a_{k}\in \mathfrak{h}$ and $m_{0}, \dots, m_{k}\in \Z_{+}$, 
\begin{eqnarray}\label{k=1-5}
\lefteqn{Y_{W}(a_{0}(-m_{0})\cdots a_{k}(-m_{k})\one, x_{1})}\nn
&&=\lim_{x_{2}\to x_{1}}\Biggl(Y_{W}(a_{0}(-m_{0})\one,x_1)
Y_{W}(a_{1}(-m_{1})\cdots a_{k}(-m_{k})\one,x_2)\nn
&&\quad\quad\quad\quad\quad-\sum_{p=1}^{k}m_{p}(a_{0}, a_{p})
{-m_{p}-1\choose m_{0}-1}(x_{1}-x_{2})^{-m_{0}-m_{p}}\cdot
\nn
&&\quad \quad \quad\quad\quad\quad\quad\quad \cdot
 Y_{W}(a_{1}(-m_{1})\cdots \widehat{a_{p}(-m_{p})}
\cdots a_{k}(-m_{k})\one,x_2)\Biggr).\nn
\end{eqnarray}
\end{cor}
\pf
Note that we can let $x_{1}=x_{2}$ in the first term of the 
right-hand side of (\ref{k=1}) and the resulting formal series is 
$$Y_{W}(a_{0}(-m_{0})\cdots a_{k}(-m_{k})\one, x_{2}).$$
Hence if we move the second term in the right-hand side of (\ref{k=1}) to the 
left-hand side, we can also let $x_{1}=x_{2}$ in the left-hand side of the 
resulting equality. Thus we obtain (\ref{k=1-5}).
\epfv

\begin{rema}
{\rm Note that in the right-hand side of (\ref{k=1-5}), we might not be able to take the 
limit (that is, let $x_{2}$ be equal to $x_{1}$) of individual terms
since we do not know whether 
the limits or substitutions exist algebraically. But the 
limit or substitution of the sum indeed exists algebraically, as is shown
in the proof of the corollary above.}
\end{rema}

We now prove the following formula for the product of two normal ordered products:

\begin{prop}
For $a_{1}, \dots, a_{k}, b_{1}, \dots, b_{l}\in \hat{\mathfrak{h}}$ and 
$m_{1}, \dots, m_{k}, n_{1}, \dots, n_{l}\in \Z_{+}$, 
\begin{eqnarray}\label{prod-normal-ordering}
\lefteqn{\nord \left(\frac{1}{(m_{1}-1)!}
\frac{\partial^{m_{1}-1}}{\partial x_{1}^{m_{1}-1}}
a_{1}(x_{1})\right)\cdots \left(\frac{1}{(m_{k}-1)!}
\frac{\partial^{m_{k}-1}}{\partial x_{k}^{m_{k}-1}}
a_{k}(x_{k})\right)\nord\cdot}\nn
&&\quad\quad\quad\cdot \nord\left(\frac{1}{(n_{1}-1)!}
\frac{\partial^{n_{1}-1}}{\partial y_{1}^{n_{1}-1}}
b_{1}(y_{1})\right)\cdots \left(\frac{1}{(n_{l}-1)!}
\frac{\partial^{n_{l}-1}}{\partial y_{l}^{n_{l}-1}}
b_{l}(y_{l})\right)\nord\nn
&&=\sum_{i=0}^{\min(k, l)}\sum_{\mbox{\scriptsize
$\begin{array}{l}k\ge p_{1}>\cdots >p_{i}\ge 1\\
0\le q_{1}<\cdots <q_{i}\le l\end{array}$}}n_{q_{1}}\cdots n_{q_{i}}(a_{p_{1}}, b_{q_{1}})
\cdots (a_{p_{i}}, b_{q_{i}})\cdot\nn
&&\quad\quad\quad\quad\cdot {-n_{q_{1}}-1\choose m_{p_{1}}-1}\cdots 
{-n_{q_{i}}-1\choose m_{p_{i}}-1}\cdot\nn
&&\quad\quad\quad\quad\cdot
(x_{p_{1}}-y_{q_{1}})^{-m_{p_{1}}-n_{q_{1}}}\cdots
(x_{p_{i}}-y_{q_{i}})^{ -m_{p_{i}}-n_{q_{i}}}\cdot\nn
&&\quad\quad\quad\quad\cdot
\nord \left(\prod_{p\ne p_{1}, \dots, p_{i}}\frac{1}{(m_{p}-1)!}
\frac{\partial^{m_{p}-1}}{\partial x_{p}^{m_{p}-1}}
a_{p}(x_{p})\right)\cdot\nn
&&\quad\quad\quad\quad\quad\quad\cdot 
\left(\prod_{q\ne q_{1}, \dots, q_{i}}\frac{1}{(n_{q}-1)!}
\frac{\partial^{n_{q}-1}}{\partial y_{q}^{n_{q}-1}}
b_{q}(y_{q})\right)\nord.
\end{eqnarray}
\end{prop}
\pf
We prove (\ref{prod-normal-ordering}) using induction on $k$. When $k=0$, (\ref{prod-normal-ordering})
holds. 

Now assume that (\ref{prod-vo-exp}) holds for $k=K$. We prove (\ref{prod-normal-ordering})
in the case $k=K+1$. For notational convenience, instead of (\ref{prod-normal-ordering}) in the case
of $k=K+1$, we prove 
(\ref{prod-normal-ordering}) with $a_{1}, \dots, a_{K+1}$ and $m_{1}, \dots, m_{K+1}$ replaced
by $a_{0}, \dots, a_{K}$ and $m_{0}, \dots, m_{K}$.
Since (\ref{prod-normal-ordering}) holds for $k=K$, we have
\begin{eqnarray}\label{prod-vo-exp-1}
\lefteqn{\left(\frac{1}{(m_{0}-1)!}
\frac{\partial^{m_{0}-1}}{\partial x_{0}^{m_{0}-1}}
a_{0}(x_{0})\right)}\nn
&&\quad\quad\quad\cdot \nord \left(\frac{1}{(m_{1}-1)!}
\frac{\partial^{m_{1}-1}}{\partial x_{1}^{m_{1}-1}}
a_{1}(x_{1})\right)\cdots \left(\frac{1}{(m_{k}-1)!}
\frac{\partial^{m_{k}-1}}{\partial x_{k}^{m_{k}-1}}
a_{k}(x_{k})\right)\nord \cdot\nn
&&\quad\quad\quad\cdot \nord\left(\frac{1}{(n_{1}-1)!}
\frac{\partial^{n_{1}-1}}{\partial y_{1}^{n_{1}-1}}
b_{1}(y_{1})\right)\cdots \left(\frac{1}{(n_{l}-1)!}
\frac{\partial^{n_{l}-1}}{\partial y_{l}^{n_{l}-1}}
b_{l}(y_{l})\right)\nord\nn
&&=\sum_{i=0}^{\min(K, l)}\sum_{\mbox{\scriptsize
$\begin{array}{l}K\ge p_{1}>\cdots >p_{i}\ge 1\\
0\le q_{1}<\cdots <q_{i}\le l\end{array}$}}n_{q_{1}}\cdots n_{q_{i}}(a_{p_{1}}, b_{q_{1}})
\cdots (a_{p_{i}}, b_{q_{i}})\cdot\nn
&&\quad\quad\quad\quad\cdot {-n_{q_{1}}-1\choose m_{p_{1}}-1}\cdots 
{-n_{q_{i}}-1\choose m_{p_{i}}-1}\cdot\nn
&&\quad\quad\quad\quad\cdot
(x_{p_{1}}-y_{q_{1}})^{-m_{p_{1}}-n_{q_{1}}}\cdots
(x_{p_{i}}-y_{q_{i}})^{ -m_{p_{i}}-n_{q_{i}}}\cdot\nn
&&\quad\quad\quad\quad\cdot \left(\frac{1}{(m_{0}-1)!}
\frac{\partial^{m_{0}-1}}{\partial x_{0}^{m_{0}-1}}
a_{0}(x_{0})\right)\cdot\nn
&&\quad\quad\quad\quad\cdot
\nord \left(\prod_{p\ne 0, p_{1}, \dots, p_{i}}\frac{1}{(m_{p}-1)!}
\frac{\partial^{m_{p}-1}}{\partial x_{p}^{m_{p}-1}}
a_{p}(x_{p})\right)\cdot\nn
&&\quad\quad\quad\quad\quad\quad\cdot 
\left(\prod_{q\ne q_{1}, \dots, q_{i}}\frac{1}{(n_{q}-1)!}
\frac{\partial^{n_{q}-1}}{\partial y_{q}^{n_{q}-1}}
b_{q}(y_{q})\right)\nord.
\end{eqnarray}
Using (\ref{normal-ordering}), the right-hand side of (\ref{prod-vo-exp-1})
is equal to 
\begin{eqnarray}\label{prod-vo-exp-2}
\lefteqn{\sum_{i=0}^{\min(K, l)}\sum_{\mbox{\scriptsize
$\begin{array}{l}K\ge p_{1}>\cdots >p_{i}\ge 1\\
0\le q_{1}<\cdots <q_{i}\le l\end{array}$}}n_{q_{1}}\cdots n_{q_{i}}(a_{p_{1}}, b_{q_{1}})
\cdots (a_{p_{i}}, b_{q_{i}})\cdot}\nn
&&\quad\quad\quad\quad\cdot {-n_{q_{1}}-1\choose m_{p_{1}}-1}\cdots 
{-n_{q_{i}}-1\choose m_{p_{i}}-1}\cdot\nn
&&\quad\quad\quad\quad\cdot
(x_{p_{1}}-y_{q_{1}})^{-m_{p_{1}}-n_{q_{1}}}\cdots
(x_{p_{i}}-y_{q_{i}})^{ -m_{p_{i}}-n_{q_{i}}}\cdot\nn
&&\quad\quad\quad\quad\cdot \nord \left(\frac{1}{(m_{0}-1)!}
\frac{\partial^{m_{0}-1}}{\partial x_{0}^{m_{0}-1}}
a_{0}(x_{0})\right)\cdot\nn
&&\quad\quad\quad\quad\quad\cdot
\left(\prod_{p\ne 0, p_{1}, \dots, p_{i}}\frac{1}{(m_{p}-1)!}
\frac{\partial^{m_{p}-1}}{\partial x_{p}^{m_{p}-1}}
a_{p}(x_{p})\right)\cdot\nn
&&\quad\quad\quad\quad\quad\quad\cdot 
\left(\prod_{q\ne q_{1}, \dots, q_{i}}\frac{1}{(n_{q}-1)!}
\frac{\partial^{n_{q}-1}}{\partial y_{q}^{n_{q}-1}}
b_{q}(y_{q})\right)\nord\nn
&&\quad + \sum_{i=0}^{\min(K, l)}\sum_{\mbox{\scriptsize
$\begin{array}{l}K\ge p_{1}>\cdots >p_{i}\ge 1\\
0\le q_{1}<\cdots <q_{i}\le l\end{array}$}}\sum_{s\ne 0, p_{1}, \dots, p_{i}}\nn
&&\quad\quad\quad\quad
m_{s} n_{q_{1}}\cdots n_{q_{i}}(a_{0}, a_{s})(a_{p_{1}}, b_{q_{1}})
\cdots (a_{p_{i}}, b_{q_{i}})\cdot\nn
&&\quad\quad\quad\quad\cdot {-m_{s}-1\choose m_{0}-1}
{-n_{q_{1}}-1\choose m_{p_{1}}-1}\cdots 
{-n_{q_{i}}-1\choose m_{p_{i}}-1}\cdot\nn
&&\quad\quad\quad\quad\cdot (x_{0}-x_{s})^{-m_{0}-m_{s}}
(x_{p_{1}}-y_{q_{1}})^{-m_{p_{1}}-n_{q_{1}}}\cdots
(x_{p_{i}}-y_{q_{i}})^{ -m_{p_{i}}-n_{q_{i}}}\cdot\nn
&&\quad\quad\quad\quad\cdot\nord
\left(\prod_{p\ne 0, s, p_{1}, \dots, p_{i}}\frac{1}{(m_{p}-1)!}
\frac{\partial^{m_{p}-1}}{\partial x_{p}^{m_{p}-1}}
a_{p}(x_{p})\right)\cdot\nn
&&\quad\quad\quad\quad\quad\quad\cdot 
\left(\prod_{q\ne q_{1}, \dots, q_{i}}\frac{1}{(n_{q}-1)!}
\frac{\partial^{n_{q}-1}}{\partial y_{q}^{n_{q}-1}}
b_{q}(y_{q})\right)\nord\nn
&&\quad + \sum_{i=0}^{\min(K, l)}\sum_{\mbox{\scriptsize
$\begin{array}{l}K\ge p_{1}>\cdots >p_{i}\ge 1\\
0\le q_{1}<\cdots <q_{i}\le l\end{array}$}}\sum_{t\ne  q_{1}, \dots, q_{i}}\nn
&&\quad\quad\quad\quad
n_{t} n_{q_{1}}\cdots n_{q_{i}}(a_{0}, b_{t})(a_{p_{1}}, b_{q_{1}})
\cdots (a_{p_{i}}, b_{q_{i}})\cdot\nn
&&\quad\quad\quad\quad\cdot {-n_{t}-1\choose m_{0}-1}
{-n_{q_{1}}-1\choose m_{p_{1}}-1}\cdots 
{-n_{q_{i}}-1\choose m_{p_{i}}-1}\cdot\nn
&&\quad\quad\quad\quad\cdot (x_{0}-y_{q})^{-m_{0}-n_{t}}
(x_{p_{1}}-y_{q_{1}})^{-m_{p_{1}}-n_{q_{1}}}\cdots
(x_{p_{i}}-y_{q_{i}})^{ -m_{p_{i}}-n_{q_{i}}}\cdot\nn
&&\quad\quad\quad\quad\cdot\nord
\left(\prod_{p\ne 0, p_{1}, \dots, p_{i}}\frac{1}{(m_{p}-1)!}
\frac{\partial^{m_{p}-1}}{\partial x_{p}^{m_{p}-1}}
a_{p}(x_{p})\right)\cdot\nn
&&\quad\quad\quad\quad\quad\quad\cdot \left(\prod_{q\ne t, 
q_{1}, \dots, q_{i}}\frac{1}{(n_{q}-1)!}
\frac{\partial^{n_{q}-1}}{\partial y_{q}^{n_{q}-1}}
b_{q}(y_{q})\right)\nord.
\end{eqnarray}

Since (\ref{prod-vo-exp}) holds for $k=K$, we also have
\begin{eqnarray}\label{prod-vo-exp-3}
\lefteqn{\nord \left(\prod_{p\ne s}\frac{1}{(m_{p}-1)!}
\frac{\partial^{m_{p}-1}}{\partial x_{p}^{m_{p}-1}}
a_{p}(x_{p})\right)\nord \nord \left(\prod_{q=1}^{l}\frac{1}{(n_{q}-1)!}
\frac{\partial^{n_{q}-1}}{\partial y_{q}^{n_{q}-1}}
b_{q}(y_{q})\right)\nord}\nn
&&=\sum_{i=0}^{\min(K, l)}\sum_{\mbox{\scriptsize
$\begin{array}{l}K\ge p_{1}>\cdots >s>\cdots >p_{i}\ge 1\\
0\le q_{1}<\cdots <q_{i}\le l\end{array}$}}n_{q_{1}}\cdots n_{q_{i}}(a_{p_{1}}, b_{q_{1}})
\cdots (a_{p_{i}}, b_{q_{i}})\cdot\nn
&&\quad\quad\quad\quad\cdot {-n_{q_{1}}-1\choose m_{p_{1}}-1}\cdots 
{-n_{q_{i}}-1\choose m_{p_{i}}-1}\cdot\nn
&&\quad\quad\quad\quad\cdot
(x_{p_{1}}-y_{q_{1}})^{-m_{p_{1}}-n_{q_{1}}}\cdots
(x_{p_{i}}-y_{q_{i}})^{ -m_{p_{i}}-n_{q_{i}}}\cdot\nn
&&\quad\quad\quad\quad\cdot
\nord \left(\prod_{p\ne s, p_{1}, \dots, p_{i}}\frac{1}{(m_{p}-1)!}
\frac{\partial^{m_{p}-1}}{\partial x_{p}^{m_{p}-1}}
a_{p}(x_{p})\right)\cdot\nn
&&\quad\quad\quad\quad\quad\quad\cdot 
\left(\prod_{q\ne q_{1}, \dots, q_{i}}\frac{1}{(n_{q}-1)!}
\frac{\partial^{n_{q}-1}}{\partial y_{q}^{n_{q}-1}}
b_{q}(y_{q})\right)\nord
\end{eqnarray}
for $s=1, \dots, k$. 

From the calculations given by (\ref{prod-vo-exp-1}), (\ref{prod-vo-exp-2}) and
(\ref{prod-vo-exp-3})
we obtain
\begin{eqnarray}\label{prod-vo-exp-4}
\lefteqn{\left(\frac{1}{(m_{0}-1)!}
\frac{\partial^{m_{0}-1}}{\partial x_{0}^{m_{0}-1}}
a_{0}(x_{0})\right)}\nn
&&\quad\quad\cdot \nord \left(\frac{1}{(m_{1}-1)!}
\frac{\partial^{m_{1}-1}}{\partial x_{1}^{m_{1}-1}}
a_{1}(x_{1})\right)\cdots \left(\frac{1}{(m_{k}-1)!}
\frac{\partial^{m_{k}-1}}{\partial x_{k}^{m_{k}-1}}
a_{k}(x_{k})\right)\nord \cdot\nn
&&\quad\quad\cdot \nord\left(\frac{1}{(n_{1}-1)!}
\frac{\partial^{n_{1}-1}}{\partial y_{1}^{n_{1}-1}}
b_{1}(y_{1})\right)\cdots \left(\frac{1}{(n_{l}-1)!}
\frac{\partial^{n_{l}-1}}{\partial y_{l}^{n_{l}-1}}
b_{l}(y_{l})\right)\nord\nn
&&\quad - \sum_{s=1}^{K}
m_{s} (a_{0}, a_{s}) {-m_{s}-1\choose m_{0}-1}  (x-x_{1})^{-m_{0}-m_{s}}\cdot\nn
&&\quad\quad\cdot
\nord \left(\prod_{p\ne s}\frac{1}{(m_{p}-1)!}
\frac{\partial^{m_{p}-1}}{\partial x_{p}^{m_{p}-1}}
a_{p}(x_{p})\right)\nord \nord \left(\prod_{q=1}^{l}\frac{1}{(n_{q}-1)!}
\frac{\partial^{n_{q}-1}}{\partial y_{q}^{n_{q}-1}}
b_{q}(y_{q})\right)\nord\nn
&&=\sum_{i=0}^{\min(K, l)}\sum_{\mbox{\scriptsize
$\begin{array}{l}K\ge p_{1}>\cdots >p_{i}\ge 1\\
0\le q_{1}<\cdots <q_{i}\le l\end{array}$}}n_{q_{1}}\cdots n_{q_{i}}(a_{p_{1}}, b_{q_{1}})
\cdots (a_{p_{i}}, b_{q_{i}})\cdot\nn
&&\quad\quad\quad\quad\cdot {-n_{q_{1}}-1\choose m_{p_{1}}-1}\cdots 
{-n_{q_{i}}-1\choose m_{p_{i}}-1}\cdot\nn
&&\quad\quad\quad\quad\cdot
(x_{p_{1}}-y_{q_{1}})^{-m_{p_{1}}-n_{q_{1}}}\cdots
(x_{p_{i}}-y_{q_{i}})^{ -m_{p_{i}}-n_{q_{i}}}\cdot\nn
&&\quad\quad\quad\quad\cdot \nord \left(\frac{1}{(m_{0}-1)!}
\frac{\partial^{m_{0}-1}}{\partial x_{0}^{m_{0}-1}}
a_{0}(x_{0})\right)\cdot\nn
&&\quad\quad\quad\quad\quad\cdot
\left(\prod_{p\ne 0, p_{1}, \dots, p_{i}}\frac{1}{(m_{p}-1)!}
\frac{\partial^{m_{p}-1}}{\partial x_{p}^{m_{p}-1}}
a_{p}(x_{p})\right)\cdot\nn
&&\quad\quad\quad\quad\quad\quad\cdot 
\left(\prod_{q\ne q_{1}, \dots, q_{i}}\frac{1}{(n_{q}-1)!}
\frac{\partial^{n_{q}-1}}{\partial y_{q}^{n_{q}-1}}
b_{q}(y_{q})\right)\nord\nn
&&\quad + \sum_{i=0}^{\min(K, l)}\sum_{\mbox{\scriptsize
$\begin{array}{l}K\ge p_{1}>\cdots >p_{i}\ge 1\\
0\le q_{1}<\cdots <q_{i}\le l\end{array}$}}\sum_{t\ne  q_{1}, \dots, q_{i}}\nn
&&\quad\quad\quad\quad
n_{t} n_{q_{1}}\cdots n_{q_{i}}(a_{0}, b_{t})(a_{p_{1}}, b_{q_{1}})
\cdots (a_{p_{i}}, b_{q_{i}})\cdot\nn
&&\quad\quad\quad\quad\cdot {-n_{t}-1\choose m_{0}-1}
{-n_{q_{1}}-1\choose m_{p_{1}}-1}\cdots 
{-n_{q_{i}}-1\choose m_{p_{i}}-1}\cdot\nn
&&\quad\quad\quad\quad\cdot (x_{0}-x_{t})^{-m_{0}-n_{t}}
(x_{p_{1}}-y_{q_{1}})^{-m_{p_{1}}-n_{q_{1}}}\cdots
(x_{p_{i}}-y_{q_{i}})^{ -m_{p_{i}}-n_{q_{i}}}\cdot\nn
&&\quad\quad\quad\quad\cdot\nord
\left(\prod_{p\ne 0, p_{1}, \dots, p_{i}}\frac{1}{(m_{p}-1)!}
\frac{\partial^{m_{p}-1}}{\partial x_{p}^{m_{p}-1}}
a_{p}(x_{p})\right)\cdot\nn
&&\quad\quad\quad\quad\quad\quad\cdot \left(\prod_{q\ne t, 
q_{1}, \dots, q_{i}}\frac{1}{(n_{q}-1)!}
\frac{\partial^{n_{q}-1}}{\partial y_{q}^{n_{q}-1}}
b_{q}(y_{q})\right)\nord\nn
&&=\sum_{i=0}^{\min(K+1, l)}\sum_{\mbox{\scriptsize
$\begin{array}{l}K\ge p_{1}>\cdots >p_{i}\ge 0\\
0\le q_{1}<\cdots <q_{i}\le l\end{array}$}}n_{q_{1}}\cdots n_{q_{i}}(a_{p_{1}}, b_{q_{1}})
\cdots (a_{p_{i}}, b_{q_{i}})\cdot\nn
&&\quad\quad\quad\quad\cdot {-n_{q_{1}}-1\choose m_{p_{1}}-1}\cdots 
{-n_{q_{i}}-1\choose m_{p_{i}}-1}\nn
&&\quad\quad\quad\quad\cdot
(x_{p_{1}}-y_{q_{1}})^{-m_{p_{1}}-n_{q_{1}}}\cdots
(x_{p_{i}}-y_{q_{i}})^{ -m_{p_{i}}-n_{q_{i}}}\cdot\nn
&&\quad\quad\quad\quad\cdot \nord \left(\frac{1}{(m_{0}-1)!}
\frac{\partial^{m_{0}-1}}{\partial x_{0}^{m_{0}-1}}
a_{0}(x_{0})\right)\cdot\nn
&&\quad\quad\quad\quad\cdot
\left(\prod_{p\ne 0, p_{1}, \dots, p_{i}}\frac{1}{(m_{p}-1)!}
\frac{\partial^{m_{p}-1}}{\partial x_{p}^{m_{p}-1}}
a_{p}(x_{p})\right)\cdot\nn
&&\quad\quad\quad\quad\quad\quad\cdot 
\left(\prod_{q\ne q_{1}, \dots, q_{i}}\frac{1}{(n_{q}-1)!}
\frac{\partial^{n_{q}-1}}{\partial x_{q}^{n_{q}-1}}
b_{q}(x_{q})\right)\nord.
\end{eqnarray}

By (\ref{normal-ordering}), the left-hand side of (\ref{prod-vo-exp-4})
is equal to 
\begin{eqnarray*}\label{prod-vo-exp-5}
\lefteqn{\nord\left(\frac{1}{(m_{0}-1)!}
\frac{\partial^{m_{0}-1}}{\partial x_{0}^{m_{0}-1}}
a_{0}(x_{0})\right)\cdot}\nn
&&\quad\quad\quad\cdot  \left(\frac{1}{(m_{1}-1)!}
\frac{\partial^{m_{1}-1}}{\partial x_{1}^{m_{1}-1}}
a_{1}(x_{1})\right)\cdots \left(\frac{1}{(m_{k}-1)!}
\frac{\partial^{m_{k}-1}}{\partial x_{k}^{m_{k}-1}}
a_{k}(x_{k})\right)\nord \cdot\nn
&&\quad\quad\quad\cdot \nord\left(\frac{1}{(n_{1}-1)!}
\frac{\partial^{n_{1}-1}}{\partial y_{1}^{n_{1}-1}}
b_{1}(y_{1})\right)\cdots \left(\frac{1}{(n_{l}-1)!}
\frac{\partial^{n_{l}-1}}{\partial y_{l}^{n_{l}-1}}
b_{l}(y_{l})\right)\nord,
\end{eqnarray*}
proving (\ref{prod-normal-ordering}) in the case $k=K+1$.
\epfv

\begin{cor}
For $a_{1}, \dots, a_{k}, b_{1}, \dots, b_{l}\in \hat{\mathfrak{h}}$ and 
$m_{1}, \dots, m_{k}, n_{1}, \dots, n_{l}\in \Z_{+}$, 
\begin{eqnarray}\label{prod-vo-exp}
\lefteqn{Y_{W}(a_{1}(-m_{1})\cdots a_{k}(-m_{k})\one, x_{1})
Y_{W}(b_{1}(-n_{1})\cdots b_{k}(-n_{l})\one, x_{2})}\nn
&&=\sum_{i=0}^{\min(k, l)}\sum_{\mbox{\scriptsize
$\begin{array}{l}k\ge p_{1}>\cdots >p_{i}\ge 1\\
0\le q_{1}<\cdots <q_{i}\le l\end{array}$}}n_{q_{1}}\cdots n_{q_{i}}(a_{p_{1}}, b_{q_{1}})
\cdots (a_{p_{i}}, b_{q_{i}})\cdot\nn
&&\quad\quad\quad\quad\cdot {-n_{q_{1}}-1\choose m_{p_{1}}-1}\cdots 
{-n_{q_{i}}-1\choose m_{p_{i}}-1}
(x_{1}-x_{2})^{-m_{p_{1}}-n_{q_{1}}-\cdots -m_{p_{i}}-n_{q_{i}}}\cdot\nn
&&\quad\quad\quad\quad\cdot
\nord \left(\prod_{p\ne p_{1}, \dots, p_{i}}\frac{1}{(m_{p}-1)!}
\frac{\partial^{m_{p}-1}}{\partial x_{1}^{m_{p}-1}}
a_{p}(x_{1})\right)\cdot\nn
&&\quad\quad\quad\quad\quad\quad\cdot \left(\prod_{q\ne q_{1}, \dots, q_{i}}\frac{1}{(n_{q}-1)!}
\frac{\partial^{n_{q}-1}}{\partial x_{2}^{n_{q}-1}}
b_{q}(x_{2})\right)\nord.
\end{eqnarray}
\end{cor}
\pf
When $x_{2}, \dots, x_{k}$ are taken to be equal to $x_{1}$ and $y_{1}, \dots, y_{l}$
are taken to be equal to $x_{2}$, the left-hand and right-hand  sides 
of (\ref{prod-normal-ordering}) exist and are equal to
the left-hand and right-hand  sides of (\ref{prod-vo-exp}). Thus (\ref{prod-vo-exp})
holds.
\epfv

The formula (\ref{prod-normal-ordering}) can be generalized to the product of 
$n$ normal ordered products. But we do not need the explicit 
formula of the coefficients.
What we need is the following:

\begin{cor}\label{general-prod-normal-ordering}
For $a_{j}\in \hat{\mathfrak{h}}$
and 
$m_{j}\in \Z_{+}$ for $j=1, \dots, n$ and $k_{0}=0, k_{1}, \dots, k_{l-1}, k_{l}=n
\in \Z$
satisfying $k_{0}=0< k_{1}<\cdots <k_{l-1}<k_{l}$, 
\begin{equation}\label{general-prod-normal-ordering-1}
\nord \prod_{j=1}^{k_{1}}\left(\frac{1}{(m_{j}-1)!}
\frac{\partial^{m_{j}-1}}{\partial x_{j}^{m_{j}-1}}
a_{j}(x_{j})\right)\nord \cdots \nord\prod_{j=k_{l-1}+1}^{n}\left(\frac{1}{(m_{j}-1)!}
\frac{\partial^{m_{j}-1}}{\partial x_{j}^{m_{j}-1}}
a_{j}(x_{j})\right)\nord
\end{equation}
is a linear combination of formal series of the form
\begin{equation}\label{general-prod-normal-ordering-2}
\prod_{(j_{1}, j_{2})\in A}(x_{j_{1}}-x_{j_{2}})^{-m_{j_{1}}-m_{j_{2}}}
\nord \prod_{j\not\in B}\left(\frac{1}{(m_{j}-1)!}
\frac{\partial^{m_{j}-1}}{\partial x_{j}^{m_{j}-1}}
a_{j}(x_{j})\right)\nord,
\end{equation}
where $A$ is a subset of 
\begin{eqnarray*}
\lefteqn{\{(j_{1}, j_{2})\;|\; \exists\; p, q\in \Z \;\mbox{\rm such that}\;
0\le p< q< l,}\nn
&&\quad\quad\quad\quad\quad\quad\quad\quad\quad\quad\quad
 k_{p}+1\le j_{1}\le k_{p+1}, \; k_{q}+1\le j_{2}\le k_{q+1} \}
\end{eqnarray*}
and $B$ is the subset of $\{1, \dots, n\}$ consisting those 
$j\in \{1, \dots, n\}$ such that 
either $(j, j')\in A$ for some $j'$ or $(j', j)\in A$ for some $j'$. 
\end{cor}
\pf
We use induction on $l$. When $l=1$, the conclusion is certainly true.
When $l=2$, (\ref{prod-normal-ordering}) gives the linear combination 
explicitly. Assume that the conclusion is 07 when $l=L$. Then
using (\ref{prod-normal-ordering}), we see that the conclusion is also true 
when $l=L+1$.
\epfv

Taking $x_{j}$ to be $x_{p+1}$ when $k_{p}+1\le j\le k_{p+1}$ in Corollary
\ref{general-prod-normal-ordering}, we obtain:

\begin{cor}\label{general-prod-vo-exp}
For $a_{j}\in \hat{\mathfrak{h}}$
and 
$m_{j}\in \Z_{+}$ for $j=1, \dots, n$ and $k_{0}=0, k_{1}, \dots, k_{l-1}, k_{l}=n
\in \Z$
satisfying $k_{0}=0< k_{1}<\cdots <k_{l-1}<k_{l}$, 
\begin{eqnarray}\label{general-prod-vo-1}
\lefteqn{Y_{W}(a_{1}(-m_{1})\cdots a_{k_{1}}(-m_{k_{1}})\one, x_{1})
\cdot}\nn
&&\quad\quad\quad\quad\quad 
\cdots Y_{W}(a_{k_{l-1}+1}(-m_{k_{l-1}+1})\cdots 
a_{n}(-m_{n})\one, x_{l})
\end{eqnarray}
is a linear combination of formal series of the form
\begin{equation}\label{general-prod-vo-2}
\prod_{(j_{1}, j_{2})\in A}(y_{j_{1}}-y_{j_{2}})^{-m_{j_{1}}-m_{j_{2}}}
\nord \prod_{j\not\in B}\left(\frac{1}{(m_{j}-1)!}
\frac{\partial^{m_{j}-1}}{\partial y_{j}^{m_{j}-1}}
a_{j}(y_{j})\right)\nord,
\end{equation}
where in the right-hand side, $y_{j}=x_{p+1}$ when $k_{p}+1\le j\le k_{p+1}$
and $A$ and $B$ are the same as in Corollary 
\ref{general-prod-normal-ordering}. \epf
\end{cor}

\section{A meromorphic open-string vertex algebra structure on $T(\hat{\mathfrak{h}}_{-})$}

In this section, we  construct a meromorphic open-string
vertex algebra structure on the left $N(\hat{\mathfrak{h}})$-module
$T(\hat{\mathfrak{h}}_{-})\otimes \C$ where we view $\C$ as a trivial left
$T(\mathfrak{h})$-module.

The left $N(\hat{\mathfrak{h}})$-module
$T(\hat{\mathfrak{h}}_{-})\otimes \C$ is 
canonically linearly isomorphic to $T(\hat{\mathfrak{h}}_{-})$.
Let $\one=1\in T(\hat{\mathfrak{h}}_{-})$. By definition, 
$T(\hat{\mathfrak{h}}_{-})$ is spanned by elements of the form
$a_{1}(-m_{1})\cdots a_{k}(-m_{k})\one$, where 
$a_{1}, \dots, a_{k}\in \hat{\mathfrak{h}}$ and $m_{1}, \dots, m_{k}\in \Z_{+}$. 

We define the weight of $a_{1}(-m_{1})\cdots a_{k}(-m_{k})\one$
for $a_{1}, \dots, a_{k}\in
\mathfrak{h}$ and $m_{1}, \dots, m_{k}\in \mathbb{Z}_{+}$ to be 
$m_{1}+\cdots + m_{k}$. Then $T(\hat\mathfrak{h}_{-})$ becomes a 
$\Z$-graded vector space. If we denote the homogeneous subspace of 
$T(\hat\mathfrak{h}_{-})$ of weight $m$ by 
$(T(\hat\mathfrak{h}_{-}))_{(m)}$, then 
$$T(\hat\mathfrak{h}_{-})=\coprod_{m\in \Z}(T(\hat\mathfrak{h}_{-}))_{(m)}.$$
The element $\one\in T(\hat\mathfrak{h}_{-})$ is 
the vacuum of $T(\hat\mathfrak{h}_{-})$.

We have a vertex operator
map 
\begin{eqnarray*}
Y_{T(\hat\mathfrak{h}_{-})}: T(\hat\mathfrak{h}_{-})&\to& 
(\mbox{\rm End}\;T(\hat\mathfrak{h}_{-}))[[x, x^{-1}]]\nn\\
v&\mapsto &Y_{T(\hat\mathfrak{h}_{-})}(u, x)=\sum_{n\in \mathbb{Z}}u_{n}x^{-n-1}
\end{eqnarray*}
where $Y_{T(\hat\mathfrak{h}_{-})}(u, x)$ is defined in (\ref{vo-map})
with $W=T(\hat\mathfrak{h}_{-})$.

We have the following main result of the present paper:

\begin{thm}\label{h-voa-0}
The triple 
$(T(\hat\mathfrak{h}_{-}), Y_{T(\hat\mathfrak{h}_{-})}, \mathbf{1})$
defined above is a meromorphic open-string vertex algebra. In the case that 
$\mathfrak{h}$ is finite dimensional, $(T(\hat\mathfrak{h}_{-}), 
Y_{T(\hat\mathfrak{h}_{-})}, \mathbf{1})$ is a grading-restricted meromorphic 
open-string vertex algebra. 
\end{thm}
\pf
The lower bound condition, 
the identity property and the creation property are easy to verify. We omit 
the proofs. It is also clear that when $\mathfrak{h}$ is finite dimensional,
the homogeneous subspaces of $T(\hat\mathfrak{h}_{-})$ are finite dimensional. 

We first prove that the 
product (\ref{product}) is absolutely convergent in the region 
$|z_{1}|>|z_{2}|>0$ to a rational function. By Corollary \ref{general-prod-vo-exp},
\begin{eqnarray}\label{prod-0}
\lefteqn{\langle v', 
Y_{T(\hat\mathfrak{h}_{-})}(a_{1}(-m_{1})\cdots a_{k_{1}}(-m_{k_{1}})\one, z_{1})
\cdot}\nn
&&\quad\quad\quad\quad
\cdots Y_{T(\hat\mathfrak{h}_{-})}(a_{k_{l-1}+1}(-m_{k_{l-1}+1})\cdots 
a_{n}(-m_{n})\one, z_{l})v\rangle
\end{eqnarray}
is a linear combination of series of the form
\begin{eqnarray}\label{prod-2}
\lefteqn{\prod_{(j_{1}, j_{2})\in A}
(y_{j_{1}}-y_{j_{2}})^{-m_{j_{1}}-m_{j_{2}}}|_{y_{j}=z_{p+1} 
\;\mbox{\scriptsize when}\;
k_{p}+1\le j\le k_{p+1}}\cdot}\nn
&&\quad\quad\cdot
\left.\left\langle v', \nord \prod_{j\not\in B}\left(\frac{1}{(m_{j}-1)!}
\frac{\partial^{m_{j}-1}}{\partial y_{j}^{m_{j}-1}}
a_{j}(y_{j})\right)\nord v\right\rangle\right|_{y_{j}=z_{p+1} 
\;\mbox{\scriptsize when}\;
k_{p}+1\le j\le k_{p+1}},\nn
\end{eqnarray}
where $A$ and $B$ are the same as in Corollary 
\ref{general-prod-normal-ordering}. Since 
$$\left.\left\langle v', \nord \prod_{j\not\in B}\left(\frac{1}{(m_{j}-1)!}
\frac{\partial^{m_{j}-1}}{\partial y_{j}^{m_{j}-1}}
a_{j}(y_{j})\right)\nord v\right\rangle\right|_{y_{j}=z_{p+1} 
\;\mbox{\scriptsize when}\;
k_{p}+1\le j\le k_{p+1}}$$
is a 
Laurent polynomial in $z_{1}, \dots, z_{l}$, 
(\ref{prod-2}) is the expansion in the 
region $|z_{1}|>\cdots >|z_{l}|>0$ of 
a rational function in $z_{1}, \dots, z_{l}$ with the only possible poles
at $z_{i}=0$ for $i=1, \dots, l$ and $z_{i}=z_{j}$ for $i\ne j$. 
Since (\ref{prod-0}) is a linear combination of series of the form 
(\ref{prod-2}), it is also the expansion in the 
region $|z_{1}|>\cdots >|z_{l}|>0$ of 
a rational function in $z_{1}, \dots, z_{l}$ with the only possible poles
at $z_{i}=0$ for $i=1, \dots, l$ and $z_{i}=z_{j}$ for $i\ne j$. 
Thus the rationality for products of vertex operators holds. 

In particular, we have rationality for products of two vertex operators. 
But in order to prove the associativity, we need an explicit 
expression of the products of two vertex operators.
By (\ref{prod-vo-exp}),
we have
\begin{eqnarray}\label{prod-1}
\lefteqn{\langle v', Y_{T(\hat\mathfrak{h}_{-})}(u_{1}, z_{1})
Y_{T(\hat\mathfrak{h}_{-})}(u_{2}, z_{2})v\rangle}\nn
&&=\sum_{i=0}^{\min(k, l)}\sum_{\mbox{\scriptsize
$\begin{array}{l}k\ge p_{1}>\cdots >p_{i}\ge 1\\
0\le q_{1}<\cdots <q_{i}\le l\end{array}$}}n_{q_{1}}\cdots n_{q_{i}}(a_{p_{1}}, b_{q_{1}})
\cdots (a_{p_{i}}, b_{q_{i}})\cdot\nn
&&\quad\quad\quad\quad\cdot {-n_{q_{1}}-1\choose m_{p_{1}}-1}\cdots 
{-n_{q_{i}}-1\choose m_{p_{i}}-1}
(z_{1}-z_{2})^{-m_{p_{1}}-n_{q_{1}}-\cdots -m_{p_{i}}-n_{q_{i}}}\cdot\nn
&&\quad\quad\quad\quad\cdot
\left\langle v', \nord \left(\prod_{p\ne p_{1}, \dots, p_{i}}\frac{1}{(m_{p}-1)!}
\frac{\partial^{m_{p}-1}}{\partial z_{1}^{m_{p}-1}}
a_{p}(z_{1})\right)\cdot\right.\nn
&&\left.\quad\quad\quad\quad\quad\quad\quad\quad\cdot \left(\prod_{q\ne q_{1}, \dots, q_{i}}\frac{1}{(n_{q}-1)!}
\frac{\partial^{n_{q}-1}}{\partial z_{2}^{n_{q}-1}}
b_{q}(z_{2})\right)\nord v\right\rangle.
\end{eqnarray}

Next we discuss iterates of two vertex operators.
From (\ref{prod-normal-ordering}), we have 
\begin{eqnarray}
\lefteqn{Y_{T(\hat\mathfrak{h}_{-})}(a_{1}(-m_{1})\cdots a_{k}(-m_{k})\one, x_{0})
b_{1}(-n_{1})\cdots b_{l}(-n_{l})\one}\nn
&&=\res_{y_{1}}\cdots\res_{y_{l}}y_{1}^{-n_{1}}\cdots y_{l}^{-n_{l}}\cdot\nn
&&\quad\quad\quad\cdot
\nord \left(\frac{1}{(m_{1}-1)!}
\frac{\partial^{m_{1}-1}}{\partial x_{0}^{m_{1}-1}}
a_{1}(x_{0})\right)\cdots \left(\frac{1}{(m_{k}-1)!}
\frac{\partial^{m_{k}-1}}{\partial x_{0}^{m_{k}-1}}
a_{k}(x_{0})\right)\nord\cdot\nn
&&\quad\quad\quad\cdot \nord
b_{1}(y_{1})\cdots 
b_{l}(y_{l})\nord\one\nn
&&=\res_{y_{1}}\cdots\res_{y_{l}}y_{1}^{-n_{1}}\cdots y_{l}^{-n_{l}}\cdot\nn
&&\quad\quad\quad\cdot\sum_{i=0}^{\min(k, l)}\sum_{\mbox{\scriptsize
$\begin{array}{l}k\ge p_{1}>\cdots >p_{i}\ge 1\\
0\le q_{1}<\cdots <q_{i}\le l\end{array}$}}(a_{p_{1}}, b_{q_{1}})
\cdots (a_{p_{i}}, b_{q_{i}})\cdot\nn
&&\quad\quad\quad\quad\cdot {-2\choose m_{p_{1}}-1}\cdots 
{-2\choose m_{p_{i}}-1}\cdot\nn
&&\quad\quad\quad\quad\cdot
(x_{0}-y_{q_{1}})^{-m_{p_{1}}-1}\cdots
(x_{0}-y_{q_{i}})^{ -m_{p_{i}}-1}\cdot\nn
&&\quad\quad\quad\quad\cdot
\nord \left(\prod_{p\ne p_{1}, \dots, p_{i}}\frac{1}{(m_{p}-1)!}
\frac{\partial^{m_{p}-1}}{\partial x_{0}^{m_{p}-1}}
a_{p}(x_{0})\right)
\left(\prod_{q\ne q_{1}, \dots, q_{i}}
b_{q}(y_{q})\right)\nord\one\nn
&&=\sum_{i=0}^{\min(k, l)}\sum_{\mbox{\scriptsize
$\begin{array}{l}k\ge p_{1}>\cdots >p_{i}\ge 1\\
0\le q_{1}<\cdots <q_{i}\le l\end{array}$}}n_{q_{1}}\cdots n_{q_{l}}(a_{p_{1}}, b_{q_{1}})
\cdots (a_{p_{i}}, b_{q_{i}})\cdot\nn
&&\quad\quad\quad\quad\cdot {-n_{q_{1}}-1\choose m_{p_{1}}-1}\cdots 
{-n_{q_{1}}-1\choose m_{p_{i}}-1}x_{0}^{-m_{p_{1}}-n_{q_{1}}-\cdots -m_{p_{i}}-n_{q_{i}}}
\cdot\nn
&&\quad\quad\quad\quad\cdot
\nord \left(\prod_{p\ne p_{1}, \dots, p_{i}}\frac{1}{(m_{p}-1)!}
\frac{\partial^{m_{p}-1}}{\partial x_{0}^{m_{p}-1}}
a_{p}^{-}(x_{0})\right)
\left(\prod_{q\ne q_{1}, \dots, q_{i}}
b_{q}(-n_{q})\right)\nord\one\nn
&&=\sum_{i=0}^{\min(k, l)}\sum_{\mbox{\scriptsize
$\begin{array}{l}k\ge p_{1}>\cdots >p_{i}\ge 1\\
0\le q_{1}<\cdots <q_{i}\le l\end{array}$}}
n_{q_{1}}\cdots n_{q_{l}}(a_{p_{1}}, b_{q_{1}})
\cdots (a_{p_{i}}, b_{q_{i}})\cdot\nn
&&\quad\quad\quad\quad\cdot {-n_{q_{1}}-1\choose m_{p_{1}}-1}\cdots 
{-n_{q_{1}}-1\choose m_{p_{i}}-1}x_{0}^{-m_{p_{1}}-n_{q_{1}}-\cdots -m_{p_{i}}-n_{q_{i}}}
\cdot\nn
&&\quad\quad\quad\quad\cdot
\left(\prod_{p\ne p_{1}, \dots, p_{i}}\sum_{s_{p}\in \Z_{+}}
{s_{p}-1\choose m_{p}-1}
a_{p}(-s_{p})x_{0}^{s_{p}-m_{p}}\right)\cdot\nn
&&\quad\quad\quad\quad\cdot
\left(\prod_{q\ne q_{1}, \dots, q_{i}}
b_{q}(-n_{q})\right)\one.\nn
\end{eqnarray}
Thus
\begin{eqnarray}\label{iterate-1}
\lefteqn{Y_{T(\hat\mathfrak{h}_{-})}(
Y_{T(\hat\mathfrak{h}_{-})}(a_{1}(-m_{1})\cdots a_{k}(-m_{k})\one, x_{0})
b_{1}(-n_{1})\cdots b_{l}(-n_{l})\one, x_{2})}\nn
&&=\sum_{i=0}^{\min(k, l)}\sum_{\mbox{\scriptsize
$\begin{array}{l}k\ge p_{1}>\cdots >p_{i}\ge 1\\
0\le q_{1}<\cdots <q_{i}\le l\end{array}$}}
n_{q_{1}}\cdots n_{q_{l}}(a_{p_{1}}, b_{q_{1}})
\cdots (a_{p_{i}}, b_{q_{i}})\cdot\nn
&&\quad\quad\quad\quad\cdot {-n_{q_{1}}-1\choose m_{p_{1}}-1}\cdots 
{-n_{q_{1}}-1\choose m_{p_{i}}-1}x_{0}^{-m_{p_{1}}-n_{q_{1}}-\cdots -m_{p_{i}}-n_{q_{i}}}
\cdot\nn
&&\quad\quad\quad\quad\cdot
\nord\left(\prod_{p\ne p_{1}, \dots, p_{i}}\sum_{s_{p}\in \Z_{+}}
{s_{p}-1\choose m_{p}-1}\cdot\right.\nn
&&\left. \quad\quad\quad\quad\quad\quad\quad\quad\quad\quad\cdot
\left(\frac{1}{(s_{p}-1)!}
\frac{\partial^{s_{p}-1}}{\partial x_{2}^{s_{p}-1}}
a_{p}(x_{2})\right)x_{0}^{s_{p}-m_{p}}\right)\cdot\nn
&&\quad\quad\quad\quad\quad\cdot
\left(\prod_{q\ne q_{1}, \dots, q_{i}}\left(\frac{1}{(n_{q}-1)!}
\frac{\partial^{n_{q}-1}}{\partial x_{2}^{n_{q}-1}}
b_{q}(x_{2})\right)\right)\nord\nn
&&=\sum_{i=0}^{\min(k, l)}\sum_{\mbox{\scriptsize
$\begin{array}{l}k\ge p_{1}>\cdots >p_{i}\ge 1\\
0\le q_{1}<\cdots <q_{i}\le l\end{array}$}}
n_{q_{1}}\cdots n_{q_{l}}(a_{p_{1}}, b_{q_{1}})
\cdots (a_{p_{i}}, b_{q_{i}})\cdot\nn
&&\quad\quad\quad\quad\cdot {-n_{q_{1}}-1\choose m_{p_{1}}-1}\cdots 
{-n_{q_{1}}-1\choose m_{p_{i}}-1}x_{0}^{-m_{p_{1}}-n_{q_{1}}-\cdots -m_{p_{i}}-n_{q_{i}}}
\cdot\nn
&&\quad\quad\quad\quad\cdot
\nord\left(\prod_{p\ne p_{1}, \dots, p_{i}}
\frac{1}{(m_{p}-1)!}
\frac{\partial^{m_{p}-1}}{\partial x_{0}^{m_{p}-1}}\cdot\right.\nn
&&\left. \quad\quad\quad\quad\quad\quad\quad\quad\quad\quad\cdot
\sum_{s_{p}\in \Z_{+}}
\left(\frac{1}{(s_{p}-1)!}
\frac{\partial^{s_{p}-1}}{\partial x_{2}^{s_{p}-1}}
a_{p}(x_{2})\right)x_{0}^{s_{p}-1}\right)\cdot\nn
&&\quad\quad\quad\quad\quad\cdot
\left(\prod_{q\ne q_{1}, \dots, q_{i}}\left(\frac{1}{(n_{q}-1)!}
\frac{\partial^{n_{q}-1}}{\partial x_{2}^{n_{q}-1}}
b_{q}(x_{2})\right)\right)\nord\nn
&&=\sum_{i=0}^{\min(k, l)}\sum_{\mbox{\scriptsize
$\begin{array}{l}k\ge p_{1}>\cdots >p_{i}\ge 1\\
0\le q_{1}<\cdots <q_{i}\le l\end{array}$}}
n_{q_{1}}\cdots n_{q_{l}}(a_{p_{1}}, b_{q_{1}})
\cdots (a_{p_{i}}, b_{q_{i}})\cdot\nn
&&\quad\quad\quad\quad\cdot {-n_{q_{1}}-1\choose m_{p_{1}}-1}\cdots 
{-n_{q_{1}}-1\choose m_{p_{i}}-1}x_{0}^{-m_{p_{1}}-n_{q_{1}}-\cdots -m_{p_{i}}-n_{q_{i}}}
\cdot\nn
&&\quad\quad\quad\quad\cdot
\nord\left(\prod_{p\ne p_{1}, \dots, p_{i}}
\frac{1}{(m_{p}-1)!}
\frac{\partial^{m_{p}-1}}{\partial x_{0}^{m_{p}-1}}
a_{p}(x_{2}+x_{0})\right)\cdot\nn
&&\quad\quad\quad\quad\quad\cdot
\left(\prod_{q\ne q_{1}, \dots, q_{i}}\left(\frac{1}{(n_{q}-1)!}
\frac{\partial^{n_{q}-1}}{\partial x_{2}^{n_{q}-1}}
b_{q}(x_{2})\right)\right)\nord.\nn
\end{eqnarray}

For $v\in T(\hat\mathfrak{h}_{-})$, 
$v'\in T(\hat\mathfrak{h}_{-})'$, from (\ref{iterate-1})
we see that when $|z_{2}|>|z_{1}-z_{2}|>0$,
the series 
\begin{eqnarray}\label{iterate-2}
\lefteqn{\langle v', Y_{T(\hat\mathfrak{h}_{-})}(
Y_{T(\hat\mathfrak{h}_{-})}(a_{1}(-m_{1})\cdots a_{k}(-m_{k})\one, z_{1}-z_{2})
b_{1}(-n_{1})\cdots b_{l}(-n_{l})\one, z_{2})v\rangle}\nn
&&=\sum_{i=0}^{\min(k, l)}\sum_{\mbox{\scriptsize
$\begin{array}{l}k\ge p_{1}>\cdots >p_{i}\ge 1\\
0\le q_{1}<\cdots <q_{i}\le l\end{array}$}}
n_{q_{1}}\cdots n_{q_{l}}(a_{p_{1}}, b_{q_{1}})
\cdots (a_{p_{i}}, b_{q_{i}})\cdot\nn
&&\quad\quad\quad\quad\cdot {-n_{q_{1}}-1\choose m_{p_{1}}-1}\cdots 
{-n_{q_{1}}-1\choose m_{p_{i}}-1}
(z_{1}-z_{2})^{-m_{p_{1}}-n_{q_{1}}-\cdots -m_{p_{i}}-n_{q_{i}}}
\cdot\nn
&&\quad\quad\quad\quad\cdot
\left\langle v', \nord\left(\prod_{p\ne p_{1}, \dots, p_{i}}
\frac{1}{(m_{p}-1)!}
\frac{\partial^{m_{p}-1}}{\partial x_{0}^{m_{p}-1}}
a_{p}(x_{2}+x_{0})\right)\cdot\right.\nn
&&\quad\quad\quad\quad\quad\cdot
\left.\left.\left(\prod_{q\ne q_{1}, \dots, q_{i}}\left(\frac{1}{(n_{q}-1)!}
\frac{\partial^{n_{q}-1}}{\partial x_{2}^{n_{q}-1}}
b_{q}(x_{2})\right)\right)\nord v\right\rangle\right|_{x_{0}=z_{1}-z_{2},\; x_{2}=z_{2}}\nn
\end{eqnarray}
is absolutely convergent to the same rational function 
to which (\ref{prod-1}) converges to. Thus rationality for iterates of two
vertex operators and associativity 
hold.

We now prove the $\mathbf{d}$-bracket property. From the definitions of the $\Z$-grading 
on $T(\hat\mathfrak{h}_{-})$, the operator $\mathbf{d}_{T(\hat\mathfrak{h}_{-})}$ 
and $a(n)$ for $a\in \mathfrak{h}$
and $n\in \Z$, we have
$$[\mathbf{d}_{T(\hat\mathfrak{h}_{-})}, a(n)]=-n a(n).$$
Then for $a\in \mathfrak{h}$,
\begin{equation}\label{d-bracket-1}
[\mathbf{d}_{T(\hat\mathfrak{h}_{-})}, a^{\pm}(x)]=a^{\pm}(x)+x\frac{d}{dx}a^{\pm}(x).
\end{equation}
For $m\in \Z_{+}$, taking $m$-th derivatives with respect to $x$ in both sides
of (\ref{d-bracket-1}), we obtain
\begin{eqnarray}\label{d-bracket-2}
\lefteqn{\left[\mathbf{d}_{T(\hat\mathfrak{h}_{-})}, 
\frac{1}{(m-1)!}\frac{d^{m-1}}{dx^{m-1}}a^{\pm}(x)\right]}\nn
&&=m\frac{1}{(m-1)!}\frac{d^{m-1}}{dx^{m-1}}a^{\pm}(x)+x\frac{d}{dx}
\left(\frac{1}{(m-1)!}
\frac{d^{m-1}}{dx^{m-1}}a^{\pm}(x)\right).
\end{eqnarray}
We also have
\begin{eqnarray}\label{d-bracket-3}
\lefteqn{\left[\mathbf{d}_{T(\hat\mathfrak{h}_{-})}, 
\frac{1}{(m-1)!}\frac{d^{m-1}}{dx^{m-1}}(a(0)x^{-1})\right]}\nn
&&=\frac{1}{(m-1)!}\frac{d^{m-1}}{dx^{m-1}}(a(0)x^{-1})+x\frac{d}{dx}
\left(\frac{1}{(m-1)!}\frac{d^{m-1}}{dx^{m-1}}(a(0)x^{-1})\right).\nn
\end{eqnarray}
Using (\ref{d-bracket-2}), (\ref{d-bracket-3}) and (\ref{vo-exp}), 
we obtain the $\mathbf{d}_{T(\hat\mathfrak{h}_{-})}$-bracket property.

We still need to prove the $D$-derivative property and the $D$-commutator formula. 
By definition,
\begin{eqnarray}\label{vo-derivative-1}
\lefteqn{\frac{d}{dx}Y_{T(\hat\mathfrak{h}_{-})}(a_{1}(-m_{1})\cdots a_{k}(-m_{k})\one, x)}\nn
&&=\frac{d}{dx}\left(\nord \left(\frac{1}{(m_{1}-1)!}
\frac{\partial^{m_{1}-1}}{\partial x^{m_{1}-1}}
a_{1}(x)\right)\cdots \left(\frac{1}{(m_{k}-1)!}
\frac{\partial^{m_{k}-1}}{\partial x^{m_{k}-1}}
a_{k}(x)\right)\nord\right)\nn
&&=\sum_{p=1}^{k}m_{p}\nord \left(\frac{1}{(m_{1}-1)!}
\frac{\partial^{m_{1}-1}}{\partial x^{m_{1}-1}}
a_{1}(x)\right)\cdot\nn
&&\quad\quad \quad\quad\quad\quad \cdots \left(\frac{1}{(m_{p-1}-1)!}
\frac{\partial^{m_{p-1}-1}}{\partial x^{m_{p-1}-1}}
a_{p-1}(x)\right) \left(\frac{1}{m_{p}!}
\frac{\partial^{m_{p}}}{\partial x^{m_{p}}}
a_{p}(x)\right)\cdot\nn
&&\quad\quad \quad\quad\quad\quad \cdot\left(\frac{1}{(m_{p+1}-1)!}
\frac{\partial^{m_{p+1}-1}}{\partial x^{m_{p+1}-1}}
a_{p+1}(x)\right)\cdot\nn
&&\quad\quad \quad\quad\quad\quad\quad\quad \cdots
\left(\frac{1}{(m_{k}-1)!}
\frac{\partial^{m_{k}-1}}{\partial x^{m_{k}-1}}
a_{k}(x)\right)\nord\nn
&&=\sum_{p=1}^{k}m_{p}Y_{T(\hat\mathfrak{h}_{-})}(a_{1}(-m_{1})\cdots 
a_{p-1}(-m_{p-1})\cdot\nn
&&\quad\quad\quad\quad\quad\quad\quad\quad\cdot 
a_{p}(-(m_{p}+1))a_{p+1}(-m_{p+1})\cdots a_{k}(-m_{k})\one, x).
\end{eqnarray}
From  (\ref{vo-derivative-1}), we obtain
\begin{eqnarray}\label{vo-derivative-2}
\lefteqn{D_{T(\hat\mathfrak{h}_{-})}a_{1}(-m_{1})\cdots a_{k}(-m_{k})\one}\nn
&&=\lim_{x\to 0}\frac{d}{dx}Y_{T(\hat\mathfrak{h}_{-})}(a_{1}(-m_{1})
\cdots a_{k}(-m_{k})\one, x)\one\nn
&&=\sum_{p=1}^{k}m_{p}\lim_{x\to 0}Y_{T(\hat\mathfrak{h}_{-})}(a_{1}(-m_{1})\cdots 
a_{p-1}(-m_{p-1})\cdot\nn
&&\quad\quad\quad\quad\quad\quad\quad\quad\cdot 
a_{p}(-(m_{p}+1))a_{p+1}(-m_{p+1})\cdots a_{k}(-m_{k})\one, x)\nn
&&=\sum_{p=1}^{k}m_{p}a_{1}(-m_{1})\cdots 
a_{p-1}(-m_{p-1})\cdot\nn
&&\quad\quad\quad\quad\quad\quad\quad\quad\cdot 
a_{p}(-(m_{p}+1))a_{p+1}(-m_{p+1})\cdots a_{k}(-m_{k})\one.
\end{eqnarray}
From (\ref{vo-derivative-1}) and (\ref{vo-derivative-2}), we obtain 
\begin{eqnarray}\label{vo-derivative-3}
\lefteqn{\frac{d}{dx}Y_{T(\hat\mathfrak{h}_{-})}(a_{1}(-m_{1})\cdots a_{k}(-m_{k})\one, x)}\nn
&&=Y_{T(\hat\mathfrak{h}_{-})}\Biggl(\sum_{p=1}^{k}m_{p}a_{1}(-m_{1})\cdots 
a_{p-1}(-m_{p-1})\cdot\nn
&&\quad\quad\quad\quad\quad\quad\quad\quad\cdot 
a_{p}(-(m_{p}+1))a_{p+1}(-m_{p+1})\cdots a_{k}(-m_{k})\one, x\Biggr)\nn
&&=Y_{T(\hat\mathfrak{h}_{-})}(Da_{1}(-m_{1})\cdots a_{k}(-m_{k})\one, x).
\end{eqnarray}

To prove 
\begin{eqnarray}\label{vo-derivative-3.1}
\lefteqn{\frac{d}{dx}Y_{T(\hat\mathfrak{h}_{-})}(a_{1}(-m_{1})\cdots a_{k}(-m_{k})\one, x)}\nn
&&=[D_{T(\hat\mathfrak{h}_{-})}, 
Y_{T(\hat\mathfrak{h}_{-})}(a_{1}(-m_{1})\cdots a_{k}(-m_{k})\one, x)]
\end{eqnarray}
for $a_{1}, \dots, a_{k}\in \mathfrak{h}$ and $m_{1}, \dots, m_{k}\in \Z_{+}$,
we use induction on $k$. When $k=0$, (\ref{vo-derivative-3.1})
holds. We also need to prove (\ref{vo-derivative-3.1}) in the case $k=1$. 
From (\ref{vo-derivative-2}), we have
\begin{equation}\label{vo-derivative-4}
[D_{T(\hat\mathfrak{h}_{-})}, a_{1}(-m_{1})]=ma_{1}(-m-1)
\end{equation}
for $a_{1}\in \mathfrak{h}$ and $m\in \Z_{+}$. 
For $b_{1}, \dots, b_{l}\in \mathfrak{h}$ and $n_{1}, \dots, n_{l}\in \Z_{+}$,
we have
\begin{eqnarray}
\lefteqn{[D_{T(\hat\mathfrak{h}_{-})}, a_{1}(m)]b_{1}(-n_{1})\cdots b_{l}(-n_{l})\one}\nn
&&=D_{T(\hat\mathfrak{h}_{-})}a_{1}(m)b_{1}(-n_{1})\cdots b_{l}(-n_{l})\one-
a_{1}(m)D_{T(\hat\mathfrak{h}_{-})}b_{1}(-n_{1})\cdots b_{l}(-n_{l})\one\nn
&&=\sum_{p=1}^{l}m(a_{1}, b_{p})\delta_{m-n_{p}, 0}Db_{1}(-n_{1})\cdots 
\widehat{b_{p}(-n_{p})}\cdots b_{l}(-n_{l})\one\nn
&&\quad -\sum_{p=1}^{l}n_{p}a_{1}(m)b_{1}(-n_{1})\cdots 
b_{p-1}(-n_{p-1})\cdot\nn
&&\quad\quad\quad\quad\quad\quad\quad\quad\cdot 
b_{p}(-(n_{p}+1))b_{p+1}(-n_{p+1})\cdots b_{l}(-n_{l})\one\nn
&&=\sum_{p=1}^{l}\sum_{q\ne p}mn_{q}(a_{1}, b_{p})\delta_{m-n_{p}, 0}\cdot\nn
&&\quad\quad\quad\quad\quad\quad\cdot b_{1}(-n_{1})\cdots 
\widehat{b_{p}(-n_{p})}\cdots b_{q}(-(n_{q}+1)\cdots b_{l}(-n_{l})\one\nn
&&\quad -\sum_{p=1}^{l}\sum_{q\ne p}mn_{p}(a_{1}, b_{q})\delta_{m-n_{q}, 0}
b_{1}(-n_{1})\cdots \widehat{b_{q}(-n_{q})}\cdots
b_{p-1}(-n_{p-1})\cdot\nn
&&\quad\quad\quad\quad\quad\quad\quad\quad\cdot 
b_{p}(-(n_{p}+1))b_{p+1}(-n_{p+1})\cdots b_{l}(-n_{l})\one\nn
&&\quad -\sum_{p=1}^{l}mn_{p}(a_{1}, b_{p})\delta_{m-n_{p}-1, 0}
b_{1}(-n_{1})\cdots 
b_{p-1}(-n_{p-1})\cdot\nn
&&\quad\quad\quad\quad\quad\quad\quad\quad\cdot 
b_{p+1}(-n_{p+1})\cdots b_{l}(-n_{l})\one\nn
&&=-m\sum_{p=1}^{l}(m-1)(a_{1}, b_{p})\delta_{m-n_{p}-1, 0}
b_{1}(-n_{1})\cdots 
b_{p-1}(-n_{p-1})\cdot\nn
&&\quad\quad\quad\quad\quad\quad\quad\quad\cdot 
b_{p+1}(-n_{p+1})\cdots b_{l}(-n_{l})\one\nn
&&=-ma_{1}(m-1)b_{1}(-n_{1})\cdots b_{l}(-n_{l})\one.
\end{eqnarray}
Thus we obtain 
\begin{equation}\label{vo-derivative-5}
[D_{T(\hat\mathfrak{h}_{-})}, a_{1}(m)]=-ma_{1}(m-1)
\end{equation}
for $a_{1}\in \mathfrak{h}$ and $m\in \Z_{+}$. The commutator formula 
(\ref{vo-derivative-4}) says that (\ref{vo-derivative-5}) holds when $m\in -\Z_{+}$.
Clearly (\ref{vo-derivative-5}) also 
holds when $m=0$. From  (\ref{vo-derivative-5}) for $m\in \Z$,
we obtain
\begin{eqnarray}\label{vo-derivative-6}
\lefteqn{[D_{T(\hat\mathfrak{h}_{-})}, Y_{T(\hat\mathfrak{h}_{-})}(a_{1}(-m_{1})\one, x)]}\nn
&&=\left[D_{T(\hat\mathfrak{h}_{-})}, \frac{1}{(m_{1}-1)!}\frac{d^{m_{1}-1}}{dx^{m_{1}-1}}a_{1}(x)\right]
=\frac{1}{(m_{1}-1)!}\frac{d^{m_{1}}}{dx^{m_{1}}}a_{1}(x)\nn
&&=\frac{d}{dx}Y_{T(\hat\mathfrak{h}_{-})}(a_{1}(-m_{1})\one, x)
\end{eqnarray}
for $a_{1}\in \mathfrak{h}$ and $m_{1}\in \Z_{+}$, proving (\ref{vo-derivative-3.1})
in the case $k=1$. 

Now assume that (\ref{vo-derivative-3.1}) holds when $k=K$. 
For $a_{0}, a_{1}, \dots, a_{k}\in \mathfrak{h}$ and $m_{0}, 
m_{1}, \dots, m_{k}\in \Z_{+}$,
from (\ref{k=1-5}) and (\ref{vo-derivative-3.1}) in the case $k=1$ and $k=K$, we obtain
\begin{eqnarray}
\lefteqn{\frac{d}{dx}Y_{T(\hat\mathfrak{h}_{-})}(a_{0}(-m_{0})
a_{1}(-m_{1})\cdots a_{k}(-m_{k})\one, x)}\nn
&&=\lim_{x_{2}\to x}\left(\frac{d}{dx}+\frac{d}{dx_{2}}\right)\cdot\nn
&&\quad \quad \quad\cdot\Biggl(Y_{T(\hat\mathfrak{h}_{-})}(a_{0}(-m_{0})\one,x)
Y_{T(\hat\mathfrak{h}_{-})}(a_{1}(-m_{1})\cdots a_{k}(-m_{k})\one,x_2)\nn
&&\quad\quad\quad\quad\quad-\sum_{p=1}^{k}m_{p}(a_{0}, a_{p})
{-m_{p}-1\choose m_{0}-1}(x-x_{2})^{-m_{0}-m_{p}}\cdot
\nn
&&\quad \quad \quad\quad\quad\quad\quad\quad \cdot
 Y_{T(\hat\mathfrak{h}_{-})}(a_{1}(-m_{1})\cdots \widehat{a_{p}(-m_{p})}
\cdots a_{k}(-m_{k})\one,x_2)\Biggr)\nn
&&=\lim_{x_{2}\to x}
\Biggl(\frac{d}{dx}Y_{T(\hat\mathfrak{h}_{-})}(a_{0}(-m_{0})\one,x)
Y_{T(\hat\mathfrak{h}_{-})}(a_{1}(-m_{1})\cdots a_{k}(-m_{k})\one,x_2)\nn
&&\quad\quad\quad\quad\quad+ Y_{T(\hat\mathfrak{h}_{-})}(a_{0}(-m_{0})\one,x)
\frac{d}{dx_{2}}Y_{T(\hat\mathfrak{h}_{-})}(a_{1}(-m_{1})\cdots a_{k}(-m_{k})\one,x_2)\nn
&&\quad\quad\quad\quad\quad-\sum_{p=1}^{k}m_{p}(a_{0}, a_{p})
{-m_{p}-1\choose m_{0}-1}(x-x_{2})^{-m_{0}-m_{p}}\cdot
\nn
&&\quad \quad \quad\quad\quad\quad\quad\quad \cdot
\frac{d}{dx_{2}}Y_{T(\hat\mathfrak{h}_{-})}(a_{1}(-m_{1})\cdots \widehat{a_{p}(-m_{p})}
\cdots a_{k}(-m_{k})\one,x_2)\Biggr)\nn
&&=\lim_{x_{2}\to x}
\Biggl([D_{T(\hat\mathfrak{h}_{-})}, Y_{T(\hat\mathfrak{h}_{-})}(a_{0}(-m_{0})\one,x)]
Y_{T(\hat\mathfrak{h}_{-})}(a_{1}(-m_{1})\cdots a_{k}(-m_{k})\one,x_2)\nn
&&\quad\quad\quad +
Y_{T(\hat\mathfrak{h}_{-})}(a_{0}(-m_{0})\one,x)
[D_{T(\hat\mathfrak{h}_{-})}, 
Y_{T(\hat\mathfrak{h}_{-})}(a_{1}(-m_{1})\cdots a_{k}(-m_{k})\one,x_2)]\nn
&&\quad\quad\quad-\sum_{p=1}^{k}m_{p}(a_{0}, a_{p})
{-m_{p}-1\choose m_{0}-1}(x-x_{2})^{-m_{0}-m_{p}}\cdot
\nn
&&\quad \quad \quad\quad\quad \cdot
[D_{T(\hat\mathfrak{h}_{-})}, 
Y_{T(\hat\mathfrak{h}_{-})}(a_{1}(-m_{1})\cdots \widehat{a_{p}(-m_{p})}
\cdots a_{k}(-m_{k})\one,x_2)]\Biggr)\nn
&&=\Biggl[D_{T(\hat\mathfrak{h}_{-})}, \lim_{x_{2}\to x}
\Biggl(Y_{T(\hat\mathfrak{h}_{-})}(a_{0}(-m_{0})\one,x)
Y_{T(\hat\mathfrak{h}_{-})}(a_{1}(-m_{1})\cdots a_{k}(-m_{k})\one,x_2)\nn
&&\quad\quad\quad\quad\quad-\sum_{p=1}^{k}m_{p}(a_{0}, a_{p})
{-m_{p}-1\choose m_{0}-1}(x-x_{2})^{-m_{0}-m_{p}}\cdot
\nn
&&\quad \quad \quad\quad\quad\quad\quad\quad \cdot
 Y_{T(\hat\mathfrak{h}_{-})}(a_{1}(-m_{1})\cdots \widehat{a_{p}(-m_{p})}
\cdots a_{k}(-m_{k})\one,x_2)\Biggr)\Biggr]\nn
&&=[D_{T(\hat\mathfrak{h}_{-})}, Y_{T(\hat\mathfrak{h}_{-})}(a_{0}(-m_{0})
a_{1}(-m_{1})\cdots a_{k}(-m_{k})\one, x)],
\end{eqnarray}
proving (\ref{vo-derivative-3.1}) in the case $k=K$.
\epfv

\begin{rema}
{\rm The symmetric algebra $S(\hat{\mathfrak{h}}_{-})$ has a natural structure of 
vertex operator algebra (see \cite{B} and \cite{FLM}, where $S(\mathfrak{h}_{-})$ is 
constructed as a subalgebra of the vertex operator algebra associated to 
a even positive definite lattice). In particular, by Remark \ref{op-va-nd-va},
it is a grading-restricted meromorphic open-string vertex algebra. 
Let $\pi: T(\hat{\mathfrak{h}}_{-}) \to S(\hat{\mathfrak{h}}_{-})$ be 
the canonical projection. Then $\pi$ is a homomorphism of meromorphic 
open-string vertex algebras from $T(\hat{\mathfrak{h}}_{-})$ to $S(\hat{\mathfrak{h}}_{-})$.
It is clear that the kernel of a homomorphism of meromorphic 
open-string vertex algebras is a subalgebra of the first meromorphic 
open-string vertex algebra and the quotient of a meromorphic 
open-string vertex algebra by a subalgebra is a meromorphic 
open-string vertex algebra. Thus we see that $S(\hat{\mathfrak{h}}_{-})$ 
as a meromorphic 
open-string vertex algebra is isomorphic to a quotient of the meromorphic 
open-string vertex algebra $T(\hat{\mathfrak{h}}_{-})$.}
\end{rema}

\section{Left modules for the meromorphic open-string vertex operator algebra 
$T(\hat{\mathfrak{h}}_{-})$}

In this section, we introduce the notion of 
left module for a meromorphic open-string vertex operator algebra. Then
we  construct a structure of a left module for the 
meromorphic open-string
vertex algebra $T(\hat{\mathfrak{h}}_{-})$ on the left $N(\hat{\mathfrak{h}})$-module
$T(\hat{\mathfrak{h}}_{-})\otimes M$ for a left $T(\mathfrak{h})$-module $M$.

\begin{defn}
{\rm Let $(V, Y_{V}, \one)$ be a meromorphic open-string vertex algebra.
A {\it module for $V$} or a {\it $V$-module} is a $\C$-graded vector space 
$W=\coprod_{n\in \C}W_{(n)}$ (graded by {\it weights}), equipped with 
a {\it vertex operator map}
\begin{eqnarray*}
Y_{W}: V&\to& ({\rm End}\; W)[[x, x^{-1}]]\nn
u&\mapsto& Y_{W}(u, x),
\end{eqnarray*}
or equivalently, 
\begin{eqnarray*}
Y_{W}: V\otimes W&\to& W[[x, x^{-1}]]\nn
u\otimes w&\mapsto& Y_{W}(u, x)w,
\end{eqnarray*}
an operator $D_{W}$ of weight $1$, satisfying the 
following conditions:
\begin{enumerate}

\item {\it Lower bound condition}:  When $\Re{(n)}$ is sufficiently negative,
$W_{(n)}=0$ 

\item The {\it identity property}:
$Y_{W}(\one,x)=1_{W}$.

\item {\it Rationality}: For $u_{1}, \dots, u_{n}, w\in W$
and $w'\in W'$, the series 
\begin{equation}\label{product-W}
\langle w', Y_{W}(u_{1}, z_1)\cdots Y_{W}(u_{n}, z_n)w\rangle
\end{equation}
converges absolutely 
when $|z_1|>\cdots >|z_n|>0$ to a rational function in $z_{1}, \dots, z_{n}$
with the only possible poles at $z_{i}=0$ for $i=1, \dots, n$ and $z_{i}=z_{j}$ 
for $i\ne j$. For $u_{1}, u_{2}, w\in W$
and $w'\in W'$, the series 
\begin{equation}\label{iterate-W}
\langle w', Y_{W}(Y_{V}(u_{1}, z_1-z_{2})u_{2}, z_2)w\rangle
\end{equation}
converges absolutely when $|z_{2}|>|z_{1}-z_{2}|>0$ to a rational function
with the only possible poles at $z_{1}=0$, $z_{2}=0$ and $z_{1}=z_{2}$. 

\item {\it Associativity}: For $u_{1}, u_{2}, w\in W$, 
$w'\in W'$, 
\begin{equation}\label{associativity-W}
\langle w', 
Y_{W}(u_{1},z_1)Y_{W}(u_{2},z_2)w\rangle
=
\langle w', 
Y_{W}(Y_{V}(u_{1},z_{1}-z_{2})u_{2},z_2)w\rangle
\end{equation}
when  $|z_{1}|>|z_{2}|>|z_{1}-z_{2}|>0$. 

\item {\it $\mathbf{d}$-bracket property}: Let $\mathbf{d}_{W}$ be the grading 
operator on $W$, that is, $\mathbf{d}_{W}w=mw$ for $m\in \R$ and
$w\in W_{(m)}$. For $u\in V$,
\begin{equation} \label{d-com-W}
[\mathbf{d}_{W}, Y_{W}(u,x)]= Y_{W}(\mathbf{d}_{V}u,x)+x\frac{d}{dx}Y_{W}(u,x).
\end{equation}

\item The {\it $D$-derivative property} and the  {\it $D$-commutator formula}: 
For $u\in V$,
\begin{eqnarray}\label{L-1-W}
\frac{d}{dx}Y_{W}(u, x)
&=&Y_{W}(D_{V}u, x) \nn
&=&[D_{W}, Y_{W}(u, x)].
\end{eqnarray}

\end{enumerate} 

A left $V$-module is said to be {\it grading restricted} if 
$\dim W_{(n)}<\infty$ for $n\in \C$. }
\end{defn}

We denote the left $V$-module just defined by $(W, Y_{W}, D_{W})$. 

\begin{rema}\label{op-va-mo-nd-va-mo}
{\rm Let $V$ be a $\Z$-graded vertex algebra such that the $\Z$-grading is lower bounded. 
By Remark \ref{op-va-nd-va}, $V$ is a meromorphic open-string vertex algebra. 
Then a $V$-module is a left module for the meromorphic open-string vertex algebra
structure.}
\end{rema}

\begin{defn}
{\rm Let $(V, Y_{V}, \one)$ be a meromorphic open-string vertex algebra
and $(W_{1}, Y_{W_{1}}, D_{W_{1}})$ and $(W_{1}, Y_{W_{1}}, D_{W_{1}})$ 
left $V$-modules. A {\it homomorphism} or {\it module map} from 
$(W_{1}, Y_{W_{1}}, D_{W_{1}})$ to $(W_{1}, Y_{W_{1}}, D_{W_{1}})$
is a linear map $f: W_{1}\to W_{2}$ such that 
\begin{eqnarray*}
f(Y_{W_{1}}(u, x)w)&=&Y_{W_{2}}(u, x)f(w),\\
f(D_{W_{1}}w)&=&D_{W_{2}}f(w)
\end{eqnarray*}
for $u\in V$ and $w\in W_{1}$. A {\it grading-preserving homomorphism} of left 
$V$-modules is a 
homomorphism
preserving the gradings. {\it Isomorphisms} or {\it equivalences} 
({\it grading-preserving isomorphisms} or {\it grading-preserving equivalence}, respectively)
are invertible homomorphisms (grading-preserving homomorphisms, respectively). 
{\it Left submodules}
({\it grading-preserving left submodules}, respectively) of a left $V$-module are 
left $V$-modules whose underlying vector spaces
are subspaces of the left $V$-module such that the embedding maps are homomorphisms 
(grading-preserving
homomorphisms, respectively).}
\end{defn}

\begin{rema}
{\rm We also have notions of  right module and bimodule for a meromorphic 
open-string vertex algebra. 
These notions and a study of these modules and left modules 
will be given in another paper on the representation theory 
of meromorphic open-string vertex algebras.}
\end{rema}

Let $M$ be a left $T(\mathfrak{h})$-module. Then we have the left 
$N(\hat\mathfrak{h})$-module 
$W=T(\hat\mathfrak{h}_{-})\otimes M$. We have a vertex operator map 
\begin{eqnarray*}
Y_{W}: T(\hat\mathfrak{h}_{-})&\to& 
(\mbox{\rm End}\; W)[[x, x^{-1}]]\nn\\
v&\mapsto &Y_{W}(u, x)=\sum_{n\in \mathbb{Z}}u_{n}x^{-n-1},
\end{eqnarray*}
where $Y_{W}(u, x)$ is defined in (\ref{vo-map}). Assume that $M$
is $\C$-graded (graded by weights) such that elements of $T(\mathfrak{h})$ preserve
the grading and the $\C$-grading is lower bounded. For example, we can 
just define $M$ to be homogeneous with an arbitrary 
complex number as the weight. 
Then this grading on $M$ together with
the grading on $T(\hat\mathfrak{h}_{-})$ gives a grading on 
$W=T(\hat\mathfrak{h}_{-})\otimes M$. Let $D_{M}$ be an operator 
on $M$ such that $D_{M}$ is of weight $1$ with respect to the grading on $M$
and commutes with the action of the elements
of $T(\mathfrak{h})$. For example, we can take $D_{M}$ to be $0$. 
We define an operator $D_{W}$ on $W$ 
by 
$$D_{W}(u\otimes w)=D_{T(\hat\mathfrak{h}_{-})}u\otimes w+u\otimes D_{M}w$$
for $u\in T(\hat\mathfrak{h}_{-})$ and $w\in M$. 
Then we have:

\begin{thm}\label{h-weak-m-0}
The triple $(W, 
Y_{W}, D_{W})$ given above is a left module for the meromorphic open-string vertex algebra 
$T(\hat\mathfrak{h}_{-})$.
\end{thm}
\pf
The  proof is in fact completely analogous to the proof of Theorem \ref{h-voa-0}. 
We omit the proof here.
\epfv

\begin{rema}
{\rm On $M$, there are actually infinitely many lower bounded $\C$-gradings.
For example, 
for any complex number, we can let the weight of every element of $M$ be
this number. 
Let $W_{1}$ and $W_{2}$ be the left $T(\hat\mathfrak{h}_{-})$-modules
obtained from the same left $T(\mathfrak{h})$-module $M$ as above and the same
$D_{M}=0$ but with 
different lower bounded $\C$-gradings. Then $W_{1}$ and $W_{2}$ are isomorphic but in
general are not grading preserving.}
\end{rema}

\begin{rema}
{\rm If $M$ is an $S(\mathfrak{h})$-module with a lower bounded $\C$-grading
(graded by weights) such that elements of $S(\hat{\mathfrak{h}}_{-})$ preserve the weights, 
then the canonical projection $\pi: T(\mathfrak{h}) \to S(\mathfrak{h})$
gives $M$ a left $T(\mathfrak{h})$-module structure with a 
lower bounded $\C$-grading
such that elements of $T(\mathfrak{h})$ preserve the weights.
Let $D_{M}$ be an operator on $M$ such that 
$D_{M}$ is of weight $1$ with respect to the grading on $M$ and 
commutes with the actions of elements of $S(\mathfrak{h})$. 
Then $S(\hat{\mathfrak{h}}_{-})\otimes M$ is a module for the 
underlying grading-restricted vertex algebra of the vertex operator algebra 
$S(\hat{\mathfrak{h}}_{-})$. The homomorphism $\pi: T(\hat{\mathfrak{h}}_{-}) 
\to S(\hat{\mathfrak{h}}_{-})$ of meromorphic open-string vertex algebras 
gives $S(\hat{\mathfrak{h}}_{-})\otimes M$ a left 
$T(\hat{\mathfrak{h}}_{-})$-module structure.
On the other hand, by Theorem \ref{h-weak-m-0},
$T(\hat{\mathfrak{h}}_{-})\otimes M$ is a left $T(\hat{\mathfrak{h}}_{-})$-module. 
Let $1_{M}$ be the identity operator on $M$. Then the map $\pi\otimes 1_{M}:
T(\hat{\mathfrak{h}}_{-})\otimes M\to S(\hat{\mathfrak{h}}_{-})\otimes M$
is a homomorphism of left $T(\hat{\mathfrak{h}}_{-})$-modules.}
\end{rema}

\begin{rema}
{\rm We also have a construction of a right $T(\hat\mathfrak{h}_{-})$-module from 
a right $T(\mathfrak{h})$-module using the construction of left 
$T(\mathfrak{h})$-modules
above. We can also construct $T(\mathfrak{h})$-bimodules. These constructions 
will be given together with the notions of right module and bimodule and 
a study of these modules and left modules
in a paper on the representation theory of 
meromorphic 
open-string vertex algebras mentioned above.}
\end{rema}

\noindent {\small \sc Department of Mathematics, Rutgers University,
110 Frelinghuysen Rd., Piscataway, NJ 08854-8019}

\vspace{1em}

\noindent {\em E-mail address}: yzhuang@math.rutgers.edu,

\end{document}